\documentclass[a4paper]{amsart}

\usepackage{amsmath,amssymb,amsthm,amscd}
\usepackage{latexsym}
\usepackage{listings}

\newtheorem{definition}{Definition}
\newtheorem{proposition}{Proposition}
\newtheorem{lemma}{Lemma}
\newtheorem{Thm}{Theorem}
\usepackage{tikz}
\usepackage{pgfplots}
\usepackage{gensymb}
\usepackage{mathtools}

\usepackage{xcolor}
\definecolor{newcolor}{rgb}{.8,.349,.1}

\renewcommand{\vec}[1]{{\bf #1}}
\DeclareMathOperator\sign{sign}

\calclayout

\begin{document}

\title[Discrete gradients in molecular dynamics]{Discrete gradients in short-range molecular dynamics simulations}

\author{Volker Grimm}
\address{Karlsruhe Institute of Technology (KIT), Institute for Applied and Numerical Mathematics, D-76131 Karlsruhe, Germany}
\email{Volker.Grimm$@$kit.edu,Tobias.Kliesch$@$kit.edu}

\author{Tobias Kliesch}

\author{G.\,R.\,W. Quispel}
\address{Department of Mathematical and Physical Sciences, La Trobe University, Victoria 3086, Australia}
\email{R.Quispel$@$latrobe.edu.au}



\date{\today}

\begin{abstract}
Discrete gradients (DG) or more exactly discrete gradient methods are time integration schemes that are custom-built to preserve first integrals or Lyapunov functions of a given ordinary differential equation (ODE).
In conservative molecular dynamics (MD) simulations, the energy of the system is constant and therefore a first integral of motion. Hence, discrete gradient methods seem to be a natural choice as an integration scheme in conservative molecular dynamics simulations.
\end{abstract}

\maketitle

\section{Introduction}
In molecular dynamics (MD) simulations, or, more generally, in particle methods, one is usually not interested in the trajectory of a single particle but in derived and/or averaged quantities that are often related to geometric properties of the underlying equations of motion.
For a meaningful simulation, the preservation of these geometric properties is important (e.g. \cite{GNI06,LeimMDbook15,McLachQui01}).
In a microcanonical or (NVE) ensemble, important macroscopic variables are the total number of particles (N), the system's volume (V) as well as the total energy (E).
Microcanonical MD simulations should therefore preserve the total energy which corresponds to a first integral of the equations of motion (e.g. \cite{SalAval14}). In the Hamiltonian formulation, the total energy or Hamiltonian $H:\mathbb{R}^d\times \mathbb{R}^d \rightarrow \mathbb{R}$, the positions $\vec{q} \in \mathbb{R}^d$ and momenta $\vec{p} \in \mathbb{R}^d$ determine the equations of motion as $\vec{q}'=\nabla_{\vec{p}}H(\vec{q},\vec{p})$, $\vec{p}'=-\nabla_\vec{q} H(\vec{q},\vec{p})$.
In addition to the preserved energy, Hamiltonian systems possess another important structural property, the symplecticity of the Hamiltonian flow (cf. \cite{LeimReich04,SanzCal94}). Unfortunately, the symplectic structure and the total energy can generally not be preserved exactly at the same time.
A class of geometric integrators that can preserve first integrals and Lyapunov functions exactly are discrete gradient methods.
They have first been considered as energy conserving schemes (e.g. \cite{Gonzalez96,GonzalezPhd96,ItohAbe88}) and they have then been generalised to arbitrary first integrals of Hamiltonian and non-Hamiltonian systems (cf. \cite{QuiTu96}) and Lyapunov functions (cf. \cite{McLachQuiTu98,StuHum96}).
Other methods that are able to preserve integrals and Lyapunov functions are projection
methods (e.g. \cite{GriQui05}), which are related to discrete gradient methods (cf. \cite{Nortonetal15}).
In addition to many applications, discrete gradient methods can be used to preserve the energy of suitably discretized variational partial differential equations (e.g.~\cite{Celledonietal12,Dahlbyetal11,Miyatakeetal17,Perseetal21,Yaguchietal12}) as well as to ensure the dissipation of gradient systems in image processing (e.g.~\cite{Grimmetal17,GrimmGPU17,Grimmetal06,Ringholmetal18}). They have been generalised to manifolds (cf. \cite{Celledonietal20,Celledonietal18}) and inspire new ideas in smooth optimization (e.g. \cite{Ehrhardtetal18,Riisetal22}) and deep learning (e.g. \cite{Celledonietal21}).

Following this introduction, discrete gradients are discussed in section~\ref{Dgigi}.
Their use for Hamiltonian systems in MD simulations is studied in section~\ref{DGmimds}. Special discrete gradients for molecular dynamics are presented in section~\ref{sec:DGmole}. The parallelisation of DG methods is discussed in section~\ref{sec:parallel}.
Finally, a brief conclusion is given in section~\ref{sec:conclusion}.
\section{Discrete gradients in geometric integration} \label{Dgigi}
In order to remind the reader of discrete gradients, we follow \cite{McLachQuiRob99}.

\begin{definition} \label{DGdef} {\bf (Gonzalez 1996)}
Let $V: \mathbb{R}^n \rightarrow \mathbb{R}$ be continuously differentiable. The function
$\overline{\nabla}V : \mathbb{R}^n \times \mathbb{R}^n \rightarrow \mathbb{R}^n$ is a {\em discrete gradient} of $V$ iff
it is continuous and
\[
 \left\{
    \begin{array}{rcl}
	\langle \overline{\nabla}V(\vec{u},\vec{u}'),(\vec{u}'-\vec{u}) \rangle &=& V(\vec{u}')-V(\vec{u}), \\
	\overline{\nabla}V(\vec{u},\vec{u}) &=& \nabla V(\vec{u}),
    \end{array}
 \right.
 \qquad \mbox{for all} \qquad \vec{u},\vec{u}' \in \mathbb{R}^n\,.
\]
The discrete gradient is {\em symmetric}, if and only if
\[
  \overline{\nabla}V(\vec{u},\vec{u}') = \overline{\nabla}V(\vec{u}',\vec{u}) \qquad \mbox{for all} \qquad \vec{u},\vec{u}' \in \mathbb{R}^n.
\]
\end{definition}
\noindent An interpretation of the following proposition, that can be found as proposition~3.2 in \cite{McLachQuiRob99}, is, that
the component of any discrete gradient in the direction  $(\vec{u}'-\vec{u}) \slash \|\vec{u}'-\vec{u}\|$ is $(V(\vec{u}')-V(\vec{u})) \slash \|\vec{u}'-\vec{u}\|$.

\begin{proposition} $\overline{\nabla}V(\vec{u},\vec{u}')$ is a discrete gradient if and only if it is continuous and
\[
\overline{\nabla}V(\vec{u},\vec{u}')=\frac{V(\vec{u}')-V(\vec{u})}{\|\vec{u}'-\vec{u}\|^2}(\vec{u}'-\vec{u})+w(\vec{u},\vec{u}'), \qquad (\vec{u} \neq \vec{u}'),
\]
 where $w(\vec{u},\vec{u}')$ is a vector-valued function such that
\[
 \left\{
   \begin{array}{l}
      \langle w(\vec{u},\vec{u}'), (\vec{u}'-\vec{u}) \rangle=0, \qquad (\vec{u} \neq \vec{u}'), \\
      \lim_{\vec{u}'\rightarrow \vec{u}}{w(\vec{u},\vec{u}')-P_{(\vec{u}'-\vec{u})^\perp}\nabla V(\vec{u})}=0,
   \end{array}
 \right.
\]
	where $P_{(\vec{u}'-\vec{u})^\perp}$ is the projection on the space perpendicular to $(\vec{u}'-\vec{u})$.
\end{proposition}
\noindent For the Euclidean inner product, $\langle \vec{u},\vec{v} \rangle = \vec{v}^T\vec{u}$, for $\vec{u},\vec{v} \in \mathbb{R}^n$, the projection can be written as the matrix
\[
  P_{(\vec{u}'-\vec{u})^\perp}=I_n-\frac{(\vec{u}'-\vec{u})}{\|\vec{u}'-\vec{u}\|^2}(\vec{u}'-\vec{u})^T\,
\]
where $I_n$ is the $n \times n$ identity matrix and $\cdot ^T$ means transpose. The following simple fact is quite useful.
\begin{lemma} \label{lem:sumdg}
Let $\overline{\nabla}V^i$ be discrete gradients for $V^i$, $i=1,\ldots,N$. Then
\[
 \overline{\nabla}V =\sum_{i=1}^N \overline{\nabla}V^i \qquad \mbox{is a discrete gradient for} \qquad V=\sum_{i=1}^N V^i\,.
\]
\end{lemma}
\noindent We omit the obvious proof of lemma~\ref{lem:sumdg}.
Examples of well-known discrete gradients are the {\em midpoint discrete gradient} or {\em Gonzalez discrete gradient} (cf. \cite{Gonzalez96})
\begin{equation}\label{gonzdg}
   \overline{V}_{MP}(\vec{u},\vec{u}') = \nabla V \left( \frac{\vec{u}+\vec{u}'}{2} \right) 
   + 
   \frac{
     V(\vec{u}')-V(\vec{u})-\left\langle \nabla V \left( \frac{\vec{u}+\vec{u}'}{2} \right), \vec{u}'-\vec{u}\right\rangle
   }{
     \|\vec{u}-\vec{u}'\|^2
   }
   (\vec{u}'-\vec{u}), \qquad (\vec{u} \neq \vec{u}')\,,
\end{equation}
the {\em mean value discrete gradient} (cf. \cite{Hartenetal83})
\begin{equation}\label{mvdg}
  \overline{V}_{MV}(\vec{u},\vec{u}')=\int_0^1 \nabla V((1-\xi)\vec{u}+\xi \vec{u}')\,d\xi, \qquad (\vec{u} \neq \vec{u}')\,,
\end{equation}
and the {\em coordinate increment discrete gradient} (cf. \cite{ItohAbe88})
\begin{equation}\label{ItohAbedg}
  \overline{V}_{CI}(\vec{u},\vec{u}')=
  \begin{pmatrix}
          \frac{V(u_1',u_2,\cdots,u_n) - V(u_1,u_2,\cdots,u_n)}{u_1'-u_1}             \\
          \frac{V(u_1',u_2',\cdots,u_n) - V(u_1',u_2,\cdots,u_n)}{u_2'-u_2}            \\
                                     \vdots                                            \\
      \frac{V(u_1',\cdots,u_{n-1}',u_n) - V(u_1',\cdots,u_{n-1},u_n)}{u_{n-1}'-u_{n-1}} \\
     \frac{V(u_1',\cdots,u_{n-1}',u_n') - V(u_1',\cdots,u_{n-1}',u_n)}{u_n'-u_n}        \\
  \end{pmatrix}\,,
\end{equation}
where $0 \slash 0$ is understood to be $\partial V \slash \partial u_i$. The midpoint discrete gradient is symmetric, i.e. $\overline{V}_{MP}(\vec{u},\vec{u}')=\overline{V}_{MP}(\vec{u}',\vec{u})$ for all $\vec{u} \neq \vec{u}'$, as is the mean value discrete gradient. The coordinate increment discrete gradient is not symmetric, but it can be symmetrised (cf. \cite{McLarenQui08})
\[
  \overline{V}_{SI}(\vec{u},\vec{u}')= \frac{1}{2} \left( \overline{V}_{CI}(\vec{u},\vec{u}') + \overline{V}_{CI}(\vec{u}',\vec{u}) \right)\,.
\]
A discrete gradient can be used as the discretised version of the gradient of a first integral or the gradient of a Lyapunov function of an ordinary differential equation (ODE). If the ODE is smooth enough and possesses a first integral or Lyapunov function $V$ then the ODE can be written as
\begin{equation}\label{lingradform}
\vec{y}'=f(\vec{y})=A(\vec{y})\nabla V(\vec{y}), 
\end{equation}
where $A(\vec{y})$ is an equally smooth antisymmetric matrix, when $V$ is a first integral, and a negative semidefinite matrix, when $V$ is a Lyapunov function (cf. \cite{McLachQuiRob99}). A {\em discrete gradient method} then reads
\begin{equation} \label{DGmethgen}
\vec{y}_{n+1}=\vec{y}_n + \tau \tilde{A}(\vec{y}_n,\vec{y}_{n+1},\tau)\overline{\nabla}V(\vec{y}_n,\vec{y}_{n+1}), \qquad n=0,1,2,\ldots\,,
\end{equation}
where $\tilde{A}(\vec{y},\vec{y},0)=A(\vec{y})$, for consistency, and $\overline{\nabla}V$ is a discrete gradient for $V$.
If the discrete gradient is symmetric and $\tilde{A}(\vec{y}_{n+1},\vec{y}_n,-\tau)=\tilde{A}(\vec{y}_n,\vec{y}_{n+1},\tau)$ holds for all possible values, then the discrete gradient method~\eqref{DGmethgen} is called time-symmetric or self-adjoint (cf. definition~3.2 in \cite{McLachQuiTu98}).
\section{DG methods in molecular dynamics (MD) simulations} \label{DGmimds}
The equations of motion in molecular dynamics can be conveniently stated in {\em Hamiltonian} form. Given an arbitrary (smooth) Hamiltonian function $H: \mathbb{R}^d \times \mathbb{R}^d \rightarrow \mathbb{R}$ on the phase space $\mathbb{R}^d \times \mathbb{R}^d$, $d \geq 1$, the corresponding Hamiltonian equations of motion are
\begin{equation} \label{genHamilneu}
\begin{array}{rcl}
	\vec{q}' &=& \phantom{-} \nabla_\vec{p} H(\vec{q},\vec{p}), \\
   \vec{p}' &=& -\nabla_\vec{q} H(\vec{q},\vec{p}).
\end{array} 
\end{equation}
The Hamiltonian corresponds to the total energy of the system that is preserved in a conservative simulation. In order to recognize the fact that equation~\eqref{genHamilneu} is a special case of equation~\eqref{lingradform}, we set $\vec{y}=(\vec{q},\vec{p})^T$ and define the matrix
\[
 J=
 \begin{bmatrix}
   0  & I_d \\
 -I_d & 0
 \end{bmatrix} \in \mathbb{R}^{2d \times 2d}, 
\]
where $I_d$ is the identity matrix of dimension $d$. With these definitions, system~\eqref{genHamilneu} reads
\begin{equation} \label{genHamil}
  \vec{y}' = J \nabla_\vec{y} H(\vec{y}) \,. 
\end{equation}
Due to
\[
\frac{d}{dt}H(\vec{y}(t))= \nabla_\vec{y} H(\vec{y}(t))^T \vec{y}'(t) = \nabla_\vec{y} H(\vec{y}(t))^T J \nabla_\vec{y} H(\vec{y}(t)) = 0\,, 
\]
the Hamiltonian $H$ (energy) is preserved along solutions of the system.
With an arbitrary discrete gradient satisfying definition~\ref{DGdef}, the corresponding discrete gradient method reads
\begin{equation} \label{DGmethgenHam}
  \vec{y}^{n+1} = \vec{y}^n + \tau J\overline {\nabla}H (\vec{y}^n,\vec{y}^{n+1})\,,
\end{equation}
with step size $\tau > 0$.
The discrete gradient method also preserves the Hamiltonian, since definition~\ref{DGdef} gives
\begin{align*}
  H(\vec{y}^{n+1})-H(\vec{y}^n) &= \overline{\nabla}H(\vec{y}^n,\vec{y}^{n+1})^T(\vec{y}^{n+1}-\vec{y}^n) \\ 
  &\stackrel{\mbox{\eqref{DGmethgenHam}}}{=} \tau \overline{\nabla}H(\vec{y}^n,\vec{y}^{n+1})^T J \overline{\nabla}H(\vec{y}^n,\vec{y}^{n+1})
  = 0\,.
\end{align*}
The method~\eqref{DGmethgenHam} is time-symmetric, whenever the chosen discrete gradient is symmetric.

\subsection{Separable Hamiltonian systems}
When the Hamiltonian is {\em separable}, i.e.
\begin{equation}\label{sepHamiltonian}
  H(\vec{q},\vec{p})=T(\vec{p})+V(\vec{q})\,, 
\end{equation}
one can apply a different discrete gradient to $T$ or $V$, respectively, in order to obtain a discrete gradient for $H(\vec{q},\vec{p})$. System \eqref{genHamil} now reads
\begin{equation} \label{sepHamil}
\begin{bmatrix}
  \vec{q} \\ \vec{p}
\end{bmatrix}'
=
\begin{bmatrix}
  0 & I \\
 -I & 0
\end{bmatrix}
\begin{bmatrix}
  \nabla_\vec{q} V \\
  \nabla_\vec{p} T
\end{bmatrix}\,. 
\end{equation}
Given two discrete gradients $\overline{\nabla}_{\vec{q}}V$ for $V$ and $\overline{\nabla}_{\vec{p}}$ $T$ for $T$, respectively, the discrete gradient method is
\begin{equation} \label{DGmethsep}
  \begin{array}{rcl}
    \vec{q}^{n+1} &=& \vec{q}^n + \tau \overline{\nabla}_{\vec{p}}T(\vec{p}^n,\vec{p}^{n+1}), \\
    \vec{p}^{n+1} &=& \vec{p}^n - \tau \overline{\nabla}_{\vec{q}}V(\vec{q}^n,\vec{q}^{n+1}).
  \end{array}
\end{equation}
Method~\eqref{DGmethsep} exactly preserves the energy, which is noted in the following lemma.
\begin{lemma} \label{energypressepHamil} For the separable Hamiltonian system \eqref{sepHamil} and two discrete gradients
$\overline{\nabla}_{\vec{q}}V$ and $\overline{\nabla}_{\vec{p}}T$, method \eqref{DGmethsep} preserves the energy exactly, i.e.
\[
  H(\vec{q}^{n+1},\vec{p}^{n+1})=H(\vec{q}^n,\vec{p}^n), \qquad n=0,1,2,\ldots\,.
\]
\end{lemma}
\begin{proof}
From definition~\ref{DGdef} followed by \eqref{DGmethsep}, we find
\begin{align*}
 H(\vec{q}^{n+1}&,\vec{p}^{n+1}) - H(\vec{q}^n,\vec{p}^n) = T(\vec{p}^{n+1})-T(\vec{q}^n)+V(\vec{q}^{n+1})-V(\vec{q}^n) \\
 &= \overline{\nabla}_{\vec{p}}T(\vec{p}^n,\vec{p}^{n+1})^T(\vec{p}^{n+1}-\vec{p}^n) + 
    \overline{\nabla}_{\vec{q}}V(\vec{q}^n,\vec{q}^{n+1})^T(\vec{q}^{n+1}-\vec{q}^n) \\
 &= -\tau \overline{\nabla}_{\vec{p}}T(\vec{p}^n,\vec{p}^{n+1})^T\overline{\nabla}_{\vec{q}}V(\vec{q}^n,\vec{q}^{n+1}) 
    +\tau \overline{\nabla}_{\vec{q}}V(\vec{q}^n,\vec{q}^{n+1})^T\overline{\nabla}_{\vec{p}}T(\vec{p}^n,\vec{p}^{n+1}) \\
 &= 0\,. 
\end{align*}
\end{proof}
\noindent The proof of lemma~\ref{energypressepHamil} also shows that
\[
 \begin{bmatrix}
   \overline{\nabla}_{\vec{p}}T \\
   \overline{\nabla}_{\vec{q}}V
 \end{bmatrix}
 =
 \overline{\nabla}H
\]
is a discrete gradient for $H$ for any choice of a discrete gradient $\overline{\nabla}_{\vec{q}}V$ and $\overline{\nabla}_{\vec{p}}T$ whenever the Hamiltonian $H$ is separable, cf. \eqref{sepHamiltonian}.
If both discrete gradients in \eqref{DGmethsep} are symmetric then the method is time-symmetric.

In molecular dynamics, the Hamiltonian often is of the even simpler form
\begin{equation} \label{sepHamilMD}
H(\vec{q},\vec{p})=\frac{1}{2}\vec{p}^TM^{-1}\vec{p} + V(\vec{q}), \qquad \mbox{i.e.} \quad T(\vec{p})=\frac{1}{2}\vec{p}^TM^{-1}\vec{p}
\end{equation}
is a quadratic function that corresponds to the kinetic energy. $M^{-1}$ is a diagonal matrix with the inverses of the masses of the corresponding particles.
For quadratic functions, any quadratic discrete gradient basically reduces to the midpoint rule.  Choosing the midpoint discrete gradient for $\overline{\nabla}_{\vec{p}}T$, one thus obtains
\begin{equation} \label{velocityDGstandardT}
  \begin{array}{rcl}
    \vec{q}^{n+1} &=& \vec{q}^n + \tau M^{-1}\frac{\vec{p}^{n+1}+\vec{p}^n}{2}, \\
    \vec{p}^{n+1} &=& \vec{p}^n - \tau \overline{\nabla}_{\vec{q}}V(\vec{q}^n,\vec{q}^{n+1}).
  \end{array}
\end{equation}
Inserting the second equation in the first leads to the system
\begin{equation} \label{velocityDGmethod}
  \begin{array}{rcl}
    \vec{q}^{n+1} &=& \vec{q}^n + \tau M^{-1}\vec{p}^n - \frac{\tau^2}{2}M^{-1}\overline{\nabla}_{\vec{q}}V(\vec{q}^n,\vec{q}^{n+1}), \\
    \vec{p}^{n+1} &=& \vec{p}^n - \tau \overline{\nabla}_{\vec{q}}V(\vec{q}^n,\vec{q}^{n+1}),
  \end{array}
\end{equation}
which will be used for the computation.
The first equation is implicit in $\vec{q}^{n+1}$ and takes some effort to solve numerically. The momenta $\vec{p}^{n+1}$ are easily computed, once the first equation has been solved.
If the first appearance of the discrete gradient were replaced by the true gradient at $\vec{q}^n$ and the second discrete gradient by the average of the true gradients at $\vec{q}^{n+1}$ and $\vec{q}^n$, respectively, we would recover the well-known Velocity-St{\"o}rmer-Verlet method. Method~\eqref{velocityDGmethod} might therefore be called Velocity-DG method.
\begin{proposition} \label{prop2}
The method~\eqref{velocityDGstandardT} or method~\eqref{velocityDGmethod}, respectively, exactly preserves the Hamiltonian~\eqref{sepHamilMD}. If the discrete gradient $\overline{\nabla}_{\vec{q}}V$ is symmetric, then the scheme is time-symmetric (reversible). If the discrete gradient $\overline{\nabla}_{\vec{q}}V$ is symmetric and sufficiently smooth then the method is of second order for sufficiently smooth $V$.
\end{proposition}
\begin{proof}
The scheme~\eqref{velocityDGmethod} exactly preserves the Hamiltonian~\eqref{sepHamilMD}, since it is just a reformulation of scheme~\eqref{DGmethsep} for this Hamiltonian and lemma~\ref{energypressepHamil} applies.
Exchanging $\tau \leftrightarrow -\tau$, $\vec{q}^{n+1} \leftrightarrow \vec{q}^n$ and $\vec{p}^{n+1} \leftrightarrow \vec{p}^n$ shows the time-symmetry for symmetric discrete gradients. Finally, definition~\ref{DGdef} and the symmetry of the scheme show second-order accuracy (cf. Theorem 8.10 in \cite{ODEI}).
\end{proof}
A more elaborate proof of second order for the midpoint discrete gradient applied to the full system~\eqref{genHamil} can be found as the proof of Theorem~8.5.4 in \cite{StuHum96}. For separable systems, this corresponds to a special case of our proposition~\ref{prop2}.

Analogous to the elimination of the velocities (momenta) in the Verlet algorithm, one may derive a two-step formulation by adding the first line of \eqref{velocityDGmethod} to the one with negative time step $-\tau$, i.e.,
\[
  \vec{q}^{n+1}-2\vec{q}^n+\vec{q}^{n-1} = -\frac{\tau^2}{2}M^{-1}
  \left(
    \overline{\nabla}_{\vec{q}}V(\vec{q}^n,\vec{q}^{n+1})+\overline{\nabla}_{\vec{q}}V(\vec{q}^n,\vec{q}^{n-1})
  \right)\,.
\]
For the solution of the implicit equation in \eqref{velocityDGmethod}, the equation is transformed to $F(\vec{u})=\vec{0}$, where
\[
F(\vec{u})=\vec{u}-\vec{q}^n-\tau M^{-1}\vec{p}^n+\frac{\tau^2}{2}M^{-1}\overline{\nabla}_{\vec{q}}V(\vec{q}^n,\vec{u})
\]
and the Newton method is applied in the vicinity of $\vec{q}^n$. This leads to the iteration
\begin{equation} \label{Newtonmeth}
    \frac{\partial F}{\partial \vec{u}}(\vec{u}^m)
    \triangle \vec{u}^m
    =
    -F(\vec{u}^m), \qquad 
    \triangle \vec{u}^m = \vec{u}^{m+1}-\vec{u}^m, \qquad m=0,1,2,\cdots
\end{equation}
where $\vec{u}^0 \approx \vec{q}^n$ and
\begin{equation} \label{NewtonmethJac}
	\frac{\partial F}{\partial \vec{u}}(\vec{u}^m) = I+\frac{\tau^2}{2}M^{-1}\frac{\partial \overline{\nabla}_{\vec{q}}V}{\partial \vec{u}} (\vec{q}^n,\vec{u}^m) \,.
\end{equation}
In order to reduce the computational work, one can also use the simplified Newton iteration
\begin{equation} \label{SimpleNewtonmeth}
    J_F(\vec{u}^m)
    \triangle \vec{u}^m
    =
    -F(\vec{u}^m), \qquad 
    \triangle \vec{u}^m = \vec{u}^{m+1}-\vec{u}^m, \qquad m=0,1,2,\cdots
\end{equation}
where $J_F(\vec{u}^m)$ is an approximation to the full Jacobian. For example, for the midpoint discrete gradient \eqref{gonzdg},
\begin{equation} \label{simpleJac}
    J_F(\vec{u}^m) = I+\frac{\tau^2}{4}M^{-1}\frac{\partial \nabla V}{\partial \vec{u}} \big(\frac{\vec{u}^m+\vec{q}^n}{2}\big)
\end{equation}
might be used. This is just the Jacobian that would occur in the implicit midpoint rule. This way, a loop over all particles to compute the potential $V(\vec{u}^m)$, the gradient $\nabla V(\vec{u}^m)$ and the norm of the difference of the positions is avoided.
The matrices $\frac{\partial F}{\partial \vec{u}}$ and  $J_F(\vec{u}^m)$ are symmetric and the linear systems in \eqref{Newtonmeth} and \eqref{SimpleNewtonmeth}, respectively, can be solved efficiently by the conjugate gradients (CG) method (cf. \cite{HeStie52}).
Newton iterations also appear in the discretisation schemes SHAKE and RATTLE that are custom-built to pose constraints during molecular dynamics simulations (cf. \cite{Andersen83,Ryckaertetal77}). The method is often also referred to as SHAKE-and-RATTLE, since they have been found to be equivalent in \cite{LeimSkeel94}.

\section{Discrete gradients for molecular dynamics} \label{sec:DGmole}
In this section, we discuss previously known as well as new discrete gradients that are designed for use in particle and molecular dynamics. 
A basic idea for discrete gradients custom-built for molecular dynamics is to mimic the standard forces by discrete gradients. That is, for a pairwise force, the discrete gradient will be built upon the two particles. For angle forces, three particles are involved and for torsion forces, four particles. Therefore, 
we consider $N$ particles with masses $m_1,\ldots,m_N$, positions $\vec{q}_i$, $i=1,\ldots,N$,  and momenta $\vec{p}_i$, $i=1,\ldots,N$. Here, $\vec{p}_i$ or $\vec{q}_i$ are three-dimensional vectors. The momenta are given as $\vec{p}_i=m_i\vec{v}_i$, where $\vec{v}_i$ are the velocities of the particles. We designate by $\vec{r}_{ij}=\vec{q}_j-\vec{q}_i$ the vector, that points from particle $i$ to particle $j$. Its length is noted as $r_{ij}=\|\vec{q}_j-\vec{q}_i\|$, where $\|\cdot \|$ designates the Euclidean norm. The connection with the vectors used in the previous sections is given as
\[
 \vec{q}=\begin{pmatrix} \vec{q}_1 \\ \vdots \\ \vec{q}_N \end{pmatrix}\,, \qquad \qquad
 \vec{p}=\begin{pmatrix} \vec{p}_1 \\ \vdots \\ \vec{p}_N \end{pmatrix}\,.	 
\]
The vector $\vec{q}$ collects the positions of the particles, the vector $\vec{p}$ the momenta, respectively.
\subsection{Discrete gradients for pairwise forces}
A discrete gradient for pairwise forces is known for quite some time. The first appearance is due to LaBudde \& Greenspan \cite{LaBuddeGreenspan74,LaBuddeGreenspan75,LaBuddeGreenspan76}. Since then, this discrete gradient has been studied by several scientists (e.g. \cite{GonzSimo96,Reich96,SimoGonz93,Simoetal92}).
If $V(r_{ij})$ is a pairwise potential for the particles $i$ and $j$ then we obtain, with the finite difference
\begin{equation} \label{dg:pw:diff}
\Delta_{r_{ij}} V(\vec{q},\vec{q}') \coloneqq \frac{V(r_{ij}')-V(r_{ij})}{r_{ij}'-r_{ij}}\,,
\end{equation}
the discrete gradient 
\[
  \overline{\nabla}V(\vec{q},\vec{q}')=
  \begin{pmatrix}
    \overline{\nabla} V_{\vec{q}_1}(\vec{q},\vec{q}') \\
    \vdots \\
    \overline{\nabla} V_{\vec{q}_N}(\vec{q},\vec{q}')
  \end{pmatrix}	  
\]	
with the non-zero components
\begin{equation} \label{dg:pw}
\overline{\nabla}_{\vec{q}_i}V(\vec{q},\vec{q}')=-\Delta_{r_{ij}} V(\vec{q},\vec{q}') \cdot
 \frac{\vec{r}_{ij}'+\vec{r}_{ij}}{r_{ij}'+r_{ij}}\,,
    \qquad
\overline{\nabla}_{\vec{q}_j}V(\vec{q},\vec{q}')=\Delta_{r_{ij}} V(\vec{q},\vec{q}') \cdot
 \frac{\vec{r}_{ij}'+\vec{r}_{ij}}{r_{ij}'+r_{ij}}\,,
\end{equation}
and $\overline{\nabla}_{\vec{q}_k} V(\vec{q},\vec{q}')=0$ for $k\neq i,j$.
\begin{Thm} \label{thm:pwdg}
The formulas \eqref{dg:pw} form a discrete gradient for pairwise potentials $V(r_{ij})$.
\end{Thm}
\begin{proof}
We have to show the properties of definition~\ref{DGdef}.
\begin{align*}
\overline{\nabla}_{\vec{q}_i}V(\vec{q},\vec{q})
  &=
	\lim_{(\tilde{\vec{q}},\vec{q}') \rightarrow (\vec{q},\vec{q})} \overline{\nabla}_{\vec{q}_i}V(\tilde{\vec{q}},\vec{q}')
	=  \lim_{(\tilde{\vec{q}},\vec{q}') \rightarrow (\vec{q},\vec{q})}
	-\Delta_{r_{ij}} V(\tilde{\vec{q}},\vec{q}') \cdot
	\frac{\vec{r}_{ij}'+\tilde{\vec{r}}_{ij}}{r_{ij}'+\tilde{r}_{ij}} \\
  &= -V'(r_{ij})\frac{\vec{r}_{ij}}{r_{ij}}=\nabla_{\vec{q}_i}V(\vec{q})\,,
\end{align*}
and, analogously,
\begin{align*}
 \overline{\nabla}_{\vec{q}_j}V(\vec{q},\vec{q})
  &=
	\lim_{(\tilde{\vec{q}},\vec{q}') \rightarrow (\vec{q},\vec{q})} \overline{\nabla}_{\vec{q}_j}V(\tilde{\vec{q}},\vec{q}')
	=  \lim_{(\tilde{\vec{q}},\vec{q}') \rightarrow (\vec{q},\vec{q})}
	\Delta_{r_{ij}} V(\tilde{\vec{q}},\vec{q}') \cdot
	\frac{\vec{r}_{ij}'+\tilde{\vec{r}}_{ij}}{r_{ij}'+\tilde{r}_{ij}} \\
  &= V'(r_{ij})\frac{\vec{r}_{ij}}{r_{ij}}=\nabla_{\vec{q}_j}V(\vec{q})\,.
\end{align*}
We also have the property
\begin{align*}
 \overline{\nabla}V(\vec{q},\vec{q}')^T(\vec{q}'-\vec{q}) 
 &= \overline{\nabla}_{\vec{q}_i}V(\vec{q},\vec{q}')^T(\vec{q}_i'-\vec{q}_i)  +
  \overline{\nabla}_{\vec{q}_j}V(\vec{q},\vec{q}')^T(\vec{q}_j'-\vec{q}_j) \\
 &= -\Delta_{r_{ij}} V \cdot
   \left(
     \frac{\vec{r}_{ij}'+\vec{r}_{ij}}{r_{ij}'+r_{ij}}
   \right)^T (\vec{q}_i'-\vec{q}_i)
 +
   \Delta_{r_{ij}} V  \cdot
   \left(
     \frac{\vec{r}_{ij}'+\vec{r}_{ij}}{r_{ij}'+r_{ij}}
    \right)^T (\vec{q}_j'-\vec{q}_j)\\
 &= \Delta_{r_{ij}} V(\vec{q},\vec{q}') \cdot
   \left(
     \frac{\vec{r}_{ij}'+\vec{r}_{ij}}{r_{ij}'+r_{ij}}
        \right)^T (\vec{r}_{ij}'-\vec{r}_{ij})
	= \Delta_{r_{ij}} V(\vec{q},\vec{q}') (r_{ij}'-r_{ij})
	= V(r_{ij}')-V(r_{ij})\\ 
  &=V(\vec{q}')-V(\vec{q})\,.
\end{align*}
\end{proof}
\noindent The following theorem is due to LaBudde and Greenspan, cf. \cite{LaBuddeGreenspan74}. The theorem and its proof might also be found as Theorem~5.1 in \cite{GNI06}. For convenience of the reader and later reference, we restate the theorem and  prove it. Note, that the important statement of the theorem is the conservation of the energy. Several schemes, also the explicit Verlet scheme, preserve the total linear momentum and the total angular momentum. The statement about the total linear momentum and the total angular momentum in the theorem are given in order to highlight that the discrete gradient methods also preserve these quantities.
\begin{Thm} \label{thm2}
 The method~\eqref{velocityDGstandardT} for a particle system of $N$ particles with Hamiltonian
\[
 H(\vec{q},\vec{p})=\frac{1}{2}\vec{p}^TM^{-1}\vec{p} + V(\vec{q})\,, \qquad V=\sum_{k=1}^M V^k\,,
\]	
where $V^k$ are pairwise potentials,
is a second-order symmetric implicit method which conserves the energy, the total linear momentum
	$\vec{P}=\sum_{i=1}^N \vec{p}_i$ and the total angular momentum $\vec{L}=\sum_{i=1}^N \vec{q}_i \times \vec{p}_i$.
\end{Thm}
\begin{proof}
The order of the method follows from proposition~\ref{prop2}.
The preservation of the energy directly follows from theorem~\ref{thm:pwdg}, lemma~\ref{lem:sumdg}, and lemma~\ref{energypressepHamil}. 	
From \eqref{dg:pw}, we have 
\[
	\sum_{i=1}^N \overline{\nabla}V^k_{\vec{q}_i} (\vec{q},\vec{q}')=0\,, \qquad \mbox{for} \qquad k=1,\ldots,M.	
\]
Hence,
\begin{equation} \label{thm2:eq1}
  \vec{P}^{n+1} = \sum_{i=1}^N \vec{p}_i^{n+1} = \sum_{i=1}^N 
    \left( 
	\vec{p}_i + \tau \sum_{k=1}^M \overline{\nabla} V^k_{\vec{q}_i}(\vec{q}^n,\vec{q}^{n+1})	
   \right)
	= \sum_{i=1}^N \vec{p}_i^n + \tau \sum_{k=1}^M \sum_{i=1}^N \overline{\nabla} V^k_{\vec{q}_i} (\vec{q}^n,\vec{q}^{n+1}) = \sum_{i=1}^N \vec{p}_i = \vec{P}^{n}.	
\end{equation}
For a fixed pairwise potential $V^k(r_{ij})$ for particles $i$ and $j$ and arbitrary $\vec{q}'$ and $\vec{q}$, we have that
\begin{align*}
  \sum_{i=1}^N (\vec{q}_i'+\vec{q}_i) \times \overline{\nabla}_{\vec{q}_i}V^k(\vec{q},\vec{q}')
  &=
   (\vec{q}_i'+\vec{q}_i)\times \overline{\nabla}_{\vec{q}_i} V^k(\vec{q},\vec{q}') +
   (\vec{q}_j'+\vec{q}_j)\times \overline{\nabla}_{\vec{q}_j} V^k(\vec{q},\vec{q}') \\ 
  &= \left[ (\vec{q}_i'+\vec{q}_i) - (\vec{q}_j'+\vec{q}_j) \right] \times \overline{\nabla}_{\vec{q}_i} V^k(\vec{q},\vec{q}') \\
  &= -\left[ (\vec{q}_i'+\vec{q}_i) - (\vec{q}_j'+\vec{q}_j) \right] \times \Delta_{r_{ij}}V^k (\vec{q},\vec{q}')\cdot \frac{\vec{r}_{ij}'+\vec{r}_{ij}}{r_{ij}'+r_{ij}} \\
  &= \frac{\Delta_{r_{ij}}V^k (\vec{q},\vec{q}')}{r_{ij}'+r_{ij}} \left[ (\vec{q}_i'+\vec{q}_i) - (\vec{q}_j'+\vec{q}_j) \right] \times \left[ (\vec{q}_i'+\vec{q}_i) - (\vec{q}_j'+\vec{q}_j) \right] 
   = \vec{0}\,.
\end{align*}
From this, for $V$ and $\vec{q}'$, $\vec{q}$ arbitrarily chosen, we find
\begin{align*}
\sum_{i=1}^N (\vec{q}_i'+\vec{q}_i) \times \overline{\nabla}_{\vec{q}_i} V(\vec{q},\vec{q}') 
   &=
    \sum_{i=1}^N (\vec{q}_i'+\vec{q}_i) \times \left( \sum_{k=1}^M \overline{\nabla}_{\vec{q}_i} V^k(\vec{q},\vec{q}') \right)	
    = \sum_{k=1}^M \left( \sum_{i=1}^N (\vec{q}_i'+\vec{q}_i) \times \overline{\nabla}_{\vec{q}_i} V^k(\vec{q},\vec{q}')\right) = \vec{0} \,.
\end{align*}
Now we turn to method~\eqref{velocityDGstandardT}. Then we have by the equation above that
\begin{equation} \label{thm2:eq2}
  \sum_{i=1}^N (\vec{q}_i^{n+1}+\vec{q}_i^n) \times (\vec{p}_i^{n+1}-\vec{p}_i^n)
  = -\tau \sum_{i=1}^N (\vec{q}_i^{n+1}+\vec{q}_i^n) \times \overline{\nabla}_{\vec{q}_i} V(\vec{q}^n,\vec{q}^{n+1}) = \vec{0}\,.
\end{equation}
Due to 
\[
  (\vec{q}_i^{n+1}-\vec{q}_i^n) \times (\vec{p}_i^{n+1}+\vec{p}_i^n) 
  = \frac{\tau}{2m_i} (\vec{p}_i^{n+1}+\vec{p}_i^n) \times (\vec{p}_i^{n+1}+\vec{p}_i^n) = \vec{0}\,,	
\]
we also have
\begin{equation} \label{thm2:eq3}
 \sum_{i=1}^N (\vec{q}_i^{n+1}-\vec{q}_i^n) \times (\vec{p}_i^{n+1}+\vec{p}_i^n) = \vec{0}\,.	
\end{equation}
Adding \eqref{thm2:eq2} and \eqref{thm2:eq3} shows the preservation of the angular momentum.
\end{proof}

\subsection{Experiment with two Lennard--Jones particles}
\begin{figure}
\begin{center}
\begin{minipage}{0.49\textwidth}
\renewcommand\lstlistingname{Code fragment}
\begin{lstlisting}[caption={Data}, escapechar=\#, label={dgabutane}]
Positions

1 6.010216  5.000000 5.000000
2 5.000000  5.000000 5.000000

Velocities

1 0.000000  1.000000 0.000000
2 0.000000 -1.000000 0.000000
\end{lstlisting}
\end{minipage}
\hspace{-2cm}
\begin{minipage}{0.49\textwidth}
\begin{center}
\hspace*{1cm} \includegraphics[height=3cm]{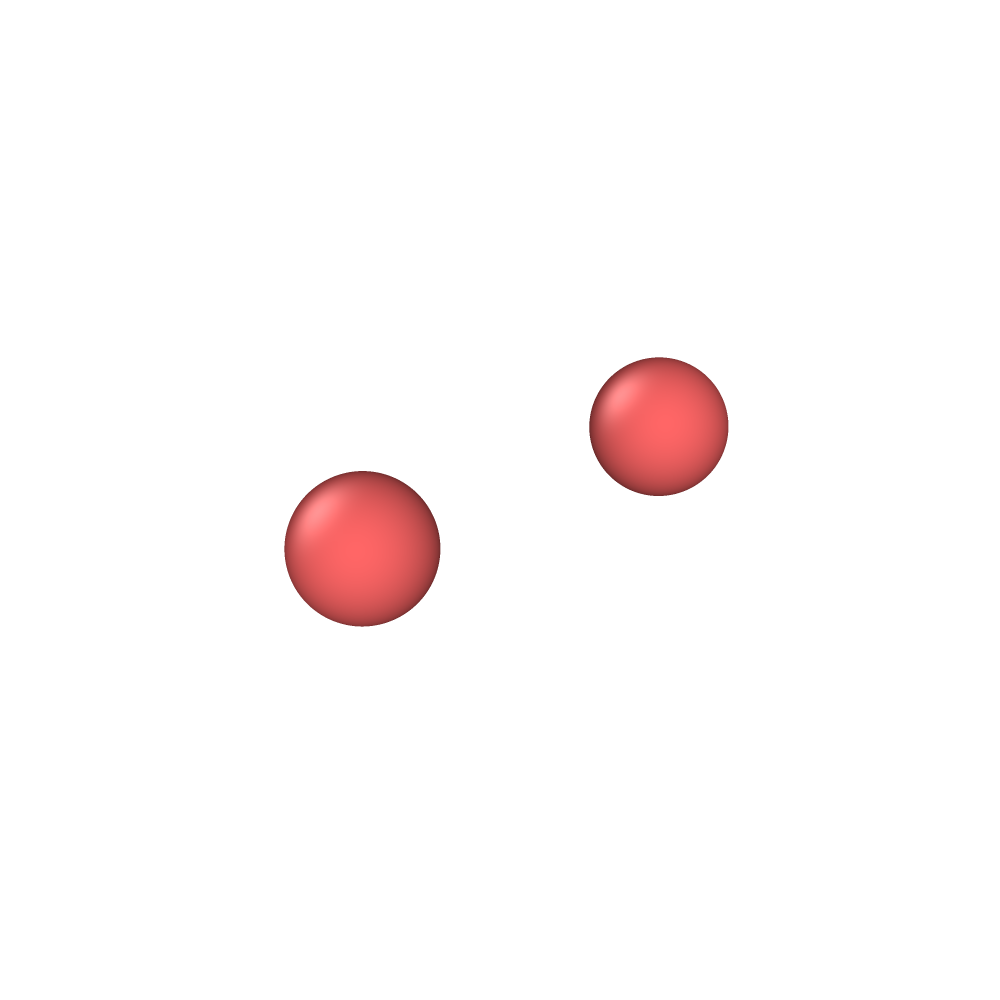}\\[4ex]
\end{center}
\end{minipage}
\end{center}
\caption{Initial conditions for the experiment with two Lennard--Jones particles: the positions and velocities are given on the left-hand side. The first row refers to the particle number, the next three rows to the three coordinates. On the right-hand side, the initial positions are given as a plot made by the open visualization tool (OVITO), cf.\cite{Stukowski10}.} \label{dgpwinitial}
\end{figure}

\begin{figure}[htp]
\begin{center}
 \begin{minipage}{0.49\textwidth}	
  \begin{tikzpicture}

\begin{axis}[
width=8cm, height=6.5cm, xmin=0.5*1e-6, xmax=2e-2, ymin=1e-10, ymax=40, xmode=log, ymode=log, legend pos=north west]
\addplot[mark=diamond*, mark size = 3pt, line width = 1pt, color=orange] coordinates {
( 1.000000000000000021e-02 , 1.830699468976740643e-02 )
( 1.000000000000000021e-03 , 1.882352545457543775e-04 )
( 1.000000000000000048e-04 , 1.882907817384800384e-06 )
( 1.000000000000000082e-05 , 1.882431550418669526e-08 )
( 9.999999999999999547e-07 , 2.162075929431185590e-10 )
};
\addplot[mark=+, mark size=3pt, line width=1pt, color=cyan] coordinates {
( 1.000000000000000021e-02 , 1.105707603241617321e+00 )
( 1.000000000000000021e-03 , 3.149699450594185129e-02 )
( 1.000000000000000048e-04 , 3.101200285280620119e-04 )
( 1.000000000000000082e-05 , 3.100953303542979620e-06 )
( 9.999999999999999547e-07 , 3.164873177202736383e-08 )
};
\addplot[mark=*, mark size=3pt, line width = 1pt, color=red!70!white] coordinates {
( 1.000000000000000021e-02 , 3.012868442463016105e+00 )
( 1.000000000000000021e-03 , 3.352693971296527575e-02 )
( 1.000000000000000048e-04 , 3.293608564546827858e-04 )
( 1.000000000000000082e-05 , 3.293220684012699715e-06 )
( 9.999999999999999547e-07 , 3.364584449723607265e-08 )
};
\addplot[color=black, thick, domain=0.000001:0.01] {100000*x*x};
\legend{DG, midpoint, Verlet, order2-line};
\end{axis}

\end{tikzpicture}
 \end{minipage}	 
 \begin{minipage}{0.49\textwidth}	
  \input{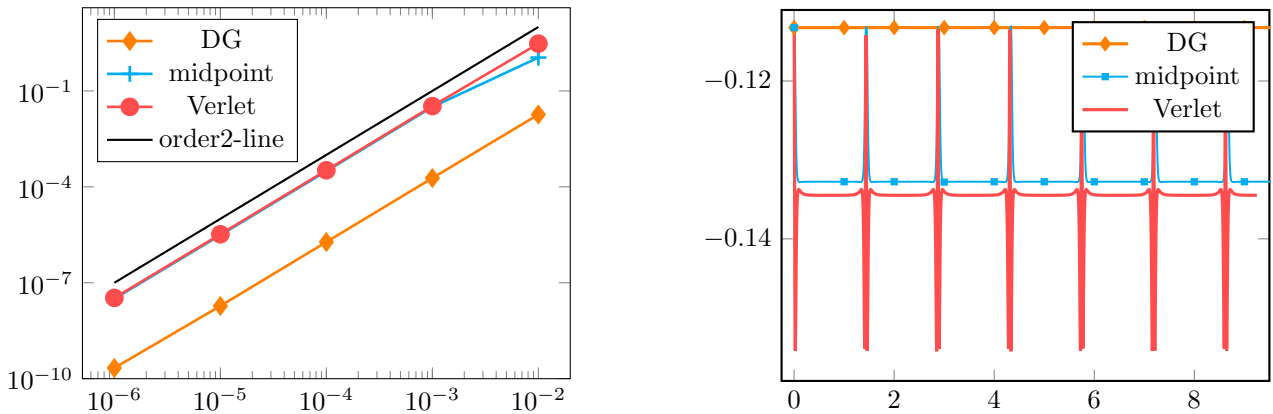}
 \end{minipage}
\caption{Results of the experiment with two Lennard--Jones particles: the error versus the time step is shown on the left-hand side for the discrete gradient (DG) method, the midpoint rule, and the Verlet scheme. On the right-hand side, the energy is shown over the time span $[0,10]$ for step-size $\tau=0.005$ for all three methods.}\label{exp:pw}
\end{center}
\end{figure}

The Lennard--Jones potential is one of the most-used models for the interaction of
neutral particles (cf. \cite{lenjones1924}). We will use the standard 12--6 potential
\begin{equation} \label{ljpotential}
 V(r_{ij})=4 \epsilon 
  \left( 
     \left( 
        \frac{\sigma}{r_{ij}}
        \right)^{12} 
     -
     \left( 
        \frac{\sigma}{r_{ij}}
     \right)^6 
  \right)
          =4 \epsilon
  \left(
     \frac{\sigma}{r_{ij}}
  \right)^6
  \cdot
  \left(
     \left(
        \frac{\sigma}{r_{ij}}
     \right)^6
     -
     1
  \right)\,.
\end{equation}
In our experiment, we use $\sigma = 1$ and $\epsilon = 5$. The initial conditions are shown in figure~\ref{dgpwinitial}. In the positions section, the atoms are numbered together with their initial location. The velocities are given in the next section, while the particle is identified by its number. The particles are attracted to each other and then repelled again. This gives a periodic movement. Due to the initial velocities, they are turning around each other, additionally.
With this set-up, we first compute the trajectories of the particles up to time $10$ with different step sizes. The error is shown in figure~\ref{exp:pw} on the left-hand side. The error is shown versus the chosen step sizes. All three methods show a second-order error behaviour. The black line indicates the second-order slope in the logarithmic plot. On the right-hand side, the calculated energies are shown for all time steps. The discrete gradient method preserves the energy correctly, while the Verlet scheme and the implicit midpoint rule show a periodic change of the energy. As expected, the total linear momentum and the total angular momentum are preserved by all three methods. 
The implicit equation in \eqref{velocityDGmethod} is solved with the Newton method, \eqref{Newtonmeth}, with the full Jacobian \eqref{NewtonmethJac}. The Jacobian is a symmetric matrix and the CG method is used to solve the linear systems. The Jacobian is not explicitly formed. Instead, the action of the Jacobian on a vector is directly computed. This way, the computational effort is comparable to computing the forces. Sometimes, this is called a matrix-free implementation.
\subsection{Discrete gradients for bond angles}
\begin{figure}[h]
\begin{center}
\begin{minipage}{0.25\textwidth}
\begin{center}
\begin{tikzpicture}[scale=1.0]
         \pgfmathsetmacro\xic{1+0.25*cos(atan(2))};
         \pgfmathsetmacro\xis{1+0.25*sin(atan(2))};
         \pgfmathsetmacro\xjc{2-0.25*sin(atan(1/2))};
         \pgfmathsetmacro\xjs{3-0.25*cos(atan(1/2))};
         \pgfmathsetmacro\xkc{3-0.25*sin(atan(1/2))};
         \pgfmathsetmacro\xks{1+0.25*cos(atan(1/2))};
         \pgfmathsetmacro\xjcr{2+0.25*sin(atan(1/2))};
         \pgfmathsetmacro\xjsr{3-0.25*cos(atan(1/2))};
         \pgfmathsetmacro\xicoff{0.05*2/sqrt(5)};
         \pgfmathsetmacro\xisoff{-0.05/sqrt(5)};
         \pgfmathsetmacro\xkcoff{0.05*2/sqrt(5)};
         \pgfmathsetmacro\xksoff{0.05/sqrt(5)};
         \draw[thick] (\xic,\xis) -- (\xjc,\xjs);
         \filldraw[white!80!black] (\xic+\xicoff,\xis+\xisoff) -- (\xjc+\xicoff,\xjs+\xisoff) -- (\xjc-\xicoff,\xjs-\xisoff) -- (\xic-\xicoff,\xis-\xisoff);
         \filldraw[white!80!black] (\xkc+\xkcoff,\xks+\xksoff) -- (\xjcr+\xkcoff,\xjsr+\xksoff) -- (\xjcr-\xkcoff,\xjsr-\xksoff) -- (\xkc-\xkcoff,\xks-\xksoff);
         \draw[thick] (1,1) circle (0.25cm);
         \node at (1,1) {$i$};
         \draw[thick] (3,1) circle (0.25cm);
         \node at (3,1) {$k$};
         \draw[thick] (2,3) circle (0.25cm);
         \node at (2,3) {$j$};
         \pgfmathsetmacro\xstart{2+1.2*cos(270-atan(1/2))};
         \pgfmathsetmacro\ystart{3+1.2*sin(270-atan(1/2))};
         \draw[<->,>=stealth,semithick] (\xstart,\ystart) arc (270-atan(1/2):270+atan(1/2):1.2cm);
	 \node at (2,1.5) {$\theta_{ijk}$};
\end{tikzpicture}
\end{center}
\end{minipage}
\hspace{-0.5cm}
\begin{minipage}{0.25\textwidth}
\begin{center}
\begin{tikzpicture}[scale=1.0]
         \pgfmathsetmacro\xic{1+0.25*cos(atan(2))};
         \pgfmathsetmacro\xis{1+0.25*sin(atan(2))};
         \pgfmathsetmacro\xjc{2-0.25*sin(atan(1/2))};
         \pgfmathsetmacro\xjs{3-0.25*cos(atan(1/2))};
         \pgfmathsetmacro\xkc{3-0.25*sin(atan(1/2))};
         \pgfmathsetmacro\xks{1+0.25*cos(atan(1/2))};
         \pgfmathsetmacro\xjcr{2+0.25*sin(atan(1/2))};
         \pgfmathsetmacro\xjsr{3-0.25*cos(atan(1/2))};
         \pgfmathsetmacro\xicoff{0.05*2/sqrt(5)};
         \pgfmathsetmacro\xisoff{-0.05/sqrt(5)};
         \pgfmathsetmacro\xkcoff{0.05*2/sqrt(5)};
         \pgfmathsetmacro\xksoff{0.05/sqrt(5)};
         \draw[thick] (\xic,\xis) -- (\xjc,\xjs);
         \filldraw[white!80!black] (\xic+\xicoff,\xis+\xisoff) -- (\xjc+\xicoff,\xjs+\xisoff) -- (\xjc-\xicoff,\xjs-\xisoff) -- (\xic-\xicoff,\xis-\xisoff);
         \filldraw[white!80!black] (\xkc+\xkcoff,\xks+\xksoff) -- (\xjcr+\xkcoff,\xjsr+\xksoff) -- (\xjcr-\xkcoff,\xjsr-\xksoff) -- (\xkc-\xkcoff,\xks-\xksoff);
         \filldraw[white!80!black] (1+0.25,1-0.05) -- (3-0.25,1-0.05) -- (3-0.25,1+0.05) -- (1+0.25,1+0.05);
         \draw[thick] (1,1) circle (0.25cm);
         \node at (1,1) {$i$};
         \draw[thick] (3,1) circle (0.25cm);
         \node at (3,1) {$k$};
         \draw[thick] (2,3) circle (0.25cm);
         \node at (2,3) {$j$};
         \pgfmathsetmacro\xstart{2+1.2*cos(270-atan(1/2))};
         \pgfmathsetmacro\ystart{3+1.2*sin(270-atan(1/2))};
         \draw[<->,>=stealth,semithick] (\xstart,\ystart) arc (270-atan(1/2):270+atan(1/2):1.2cm);
	 \node at (2,1.5) {$\theta_{ijk}$};
\end{tikzpicture}
\end{center}
\end{minipage}
\hfill
\begin{minipage}{0.49\textwidth}
\[
\hspace{-2cm}
	\theta_{ijk} =
   \arccos
   \left(
     \frac{
        \langle \vec{r}_{ji},\vec{r}_{jk}\rangle
     }{
             \|\vec{r}_{ji}\|\cdot \|\vec{r}_{jk}\|
     }
        \right) 
 =
   \arccos
   \left(
     \frac{
             r_{ji}^2+r_{jk}^2-r_{ik}^2
     }{
             2r_{ji} \cdot r_{jk}
     }
   \right)
\]
\end{minipage}
\end{center}
\caption{Sketch of bond angle: From left to right, first the bonds and the bond angle between the atoms $i$,$j$,$k$ is shown. The second sketch shows the distances between the atoms $i$,$j$,$k$ used in the discrete gradient. On the right-hand side, the standard term of the bond angle is followed by the term based on distances, solely.} \label{angleplotpic}
\end{figure}
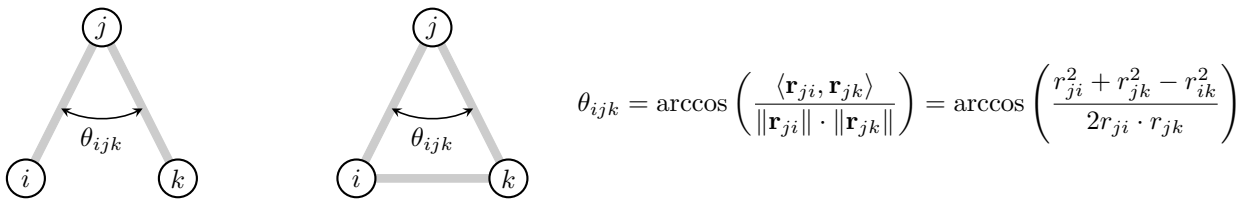

We discuss some discrete gradients for bond angles. Besides general discrete gradients, e.g. \eqref{gonzdg}, \eqref{mvdg}, and \eqref{ItohAbedg}, restricted to the bond angles, a symmetric discrete gradient for the bond angles has recently been proposed in \cite{SchieblRomero21}. We discuss a slight generalisation of this discrete gradient. 
The standard angle potential, $V(\theta)$, is assumed to depend smoothly on the angle $\theta$, cf. figure~\ref{angleplotpic}.
The angle $\theta=\theta_{ijk}$ can be expressed in terms of the distances $r_{ji}$, $r_{jk}$, $r_{ij}$ between the three atoms $i$, $j$, and $k$ as given on the right-hand side of figure~\ref{angleplotpic}. This is due to the following lemma.
\begin{lemma}
We have the following representation of the scalar product
\[
 \langle \vec{r}_{ji}, \vec{r}_{jk} \rangle 
   = \frac{1}{2}
   \left(
     r_{ji}^2+r_{jk}^2-r_{ik}^2
   \right)\,.
\]
\end{lemma}
\begin{proof}
The proof is a simple calculation.
\begin{align*}
r_{ji}^2 + r_{jk}^2 - r_{ik}^2 &=
 \langle \vec{q}_i - \vec{q}_j, \vec{q}_i - \vec{q}_j \rangle
 + \langle \vec{q}_k - \vec{q}_j, \vec{q}_k - \vec{q}_j \rangle
 - \langle \vec{q}_k - \vec{q}_i, \vec{q}_k - \vec{q}_i \rangle \\
 &=
 2 \left(
  \langle \vec{q}_j, \vec{q}_j \rangle
  - \langle \vec{q}_i, \vec{q}_j \rangle
  - \langle \vec{q}_k, \vec{q}_j \rangle
  + \langle \vec{q}_i, \vec{q}_k \rangle
 \right)
  =
  2 \cdot \langle \vec{q}_i-\vec{q}_j, \vec{q}_k-\vec{q}_j  \rangle 
  = 2 \cdot \langle \vec{r}_{ji}, \vec{r}_{jk} \rangle\,.
\end{align*}
\end{proof}
\noindent So the angle potential only depends on three distances $V(r_{ji},r_{jk},r_{ik})$. We first give a non-symmetric discrete gradient that is similar to the Itoh--Abe gradient. In order to write the discrete gradient down in a concise way, we use the following finite differences with respect to the distances of the particles
\begin{align*}
\Delta_{r_{ji}}^\ell V(\vec{q},\vec{q}')
&\coloneqq
\frac{
   V(r_{ji}',r_{jk},r_{ik})-V(r_{ji},r_{jk},r_{ik})
  }{
   r_{ji}'-r_{ji}
  }\,,\\
\Delta_{r_{jk}}^\ell V(\vec{q},\vec{q}')
&\coloneqq
\frac{
   V(r_{ji}',r_{jk}',r_{ik})-V(r_{ji}',r_{jk},r_{ik})
  }{
   r_{jk}'-r_{jk}
  }\,, \\
\Delta_{r_{ik}}^\ell V(\vec{q},\vec{q}')
&\coloneqq
 \frac{
   V(r_{ji}',r_{jk}',r_{ik}')-V(r_{ji}',r_{jk}',r_{ik})
  }{
   r_{ik}'-r_{ik}
  }\,.
\end{align*}
This can also be seen as the Itoh--Abe (coordinate increment) discrete gradient for $V(r_{ji},r_{jk},r_{ik})$. With these finite differences, our first discrete gradient $\overline{\nabla}^\ell V$ has the nonzero components
\begin{align}
\overline{\nabla}_{\vec{q}_i}^\ell V(\vec{q},\vec{q}') &=
  \phantom{-} \Delta_{r_{ji}}^\ell V(\vec{q},\vec{q}')
  \cdot
  \frac{\vec{r}_{ji}'+\vec{r}_{ji}}{r_{ji}'+r_{ji}}
  -
  \Delta_{r_{ik}}^\ell V(\vec{q},\vec{q}')
  \cdot
  \frac{\vec{r}_{ik}'+\vec{r}_{ik}}{r_{ik}'+r_{ik}}\,, \nonumber \\[1ex]
\overline{\nabla}_{\vec{q}_j}^\ell V(\vec{q},\vec{q}') &=
  -\Delta_{r_{ji}}^\ell V(\vec{q},\vec{q}')
  \cdot
  \frac{\vec{r}_{ji}'+\vec{r}_{ji}}{r_{ji}'+r_{ji}}
  -
  \Delta_{r_{jk}}^\ell V(\vec{q},\vec{q}')
  \cdot
  \frac{\vec{r}_{jk}'+\vec{r}_{jk}}{r_{jk}'+r_{jk}}\,,\label{dganglel} \\[1ex]
\overline{\nabla}_{\vec{q}_k}^\ell V(\vec{q},\vec{q}') &=
  \phantom{-} \Delta_{r_{jk}}^\ell V(\vec{q},\vec{q}')
  \cdot
  \frac{\vec{r}_{jk}'+\vec{r}_{jk}}{r_{jk}'+r_{jk}}
  +
  \Delta_{r_{ik}}^\ell V(\vec{q},\vec{q}')
  \cdot
  \frac{\vec{r}_{ik}'+\vec{r}_{ik}}{r_{ik}'+r_{ik}}\,,\nonumber 
\end{align}
where all other are zero.
This discrete gradient is not symmetric, since the finite differences are not symmetric, e.g.
\[
\Delta_{r_{ji}}^\ell V(\vec{q},\vec{q}') \neq \Delta_{r_{ji}}^\ell V(\vec{q}',\vec{q}). 
\]
Hence we can define another discrete gradient $\overline{\nabla}^r V$ by using coordinate decrease. The finite differences are now
\begin{align*}
\Delta_{r_{ji}}^r V(\vec{q},\vec{q}') \coloneqq \Delta_{r_{ji}}^\ell V(\vec{q}',\vec{q})
&=
\frac{
  V(r_{ji}',r_{jk}',r_{ik}')-V(r_{ji},r_{jk}',r_{ik}')
  }{
   r_{ji}'-r_{ji}
  }\, , \\
\Delta_{r_{jk}}^r V(\vec{q},\vec{q}') \coloneqq \Delta_{r_{jk}}^\ell V(\vec{q}',\vec{q})
&=
\frac{
  V(r_{ji},r_{jk}',r_{ik}')-V(r_{ji},r_{jk},r_{ik}')
  }{
   r_{jk}'-r_{jk}
  }\, , \\
\Delta_{r_{ik}}^r V(\vec{q},\vec{q}') \coloneqq \Delta_{r_{ik}}^\ell V(\vec{q}',\vec{q})
&=
 \frac{
   V(r_{ji},r_{jk},r_{ik}')-V(r_{ji},r_{jk},r_{ik})
  }{
   r_{ik}'-r_{ik}
  }\, . 
\end{align*}
With these finite differences, we obtain another discrete gradient that is also not symmetric.
Since symmetric discrete gradients lead to symmetric time-integration methods, which are of second order, one might prefer symmetric discrete gradients. In the same way as the Itoh--Abe discrete gradient, we can symmetrise our gradients here. This leads to the discrete gradient with the (symmetric) finite differences
\begin{align*}
 \Delta_v^s V(\vec{q},\vec{q}')
   \coloneqq \frac{1}{2}
    \left(
        \Delta_v^\ell V(\vec{q},\vec{q}')+\Delta_v^r V(\vec{q},\vec{q}')
	\right), \qquad v \in \{ r_{ji}, r_{jk}, r_{ik} \}.
\end{align*}
And the non-zero coefficients of the symmetric discrete gradient $\overline{\nabla}^s V$ read
\begin{align}
\overline{\nabla}_{\vec{q}_i}^sV(\vec{q},\vec{q}') &=
  \phantom{-} \Delta_{r_{ji}}^s V(\vec{q},\vec{q}')
  \cdot
  \frac{\vec{r}_{ji}'+\vec{r}_{ji}}{r_{ji}'+r_{ji}}
  -
  \Delta_{r_{ik}}^s V(\vec{q},\vec{q}')
  \cdot
  \frac{\vec{r}_{ik}'+\vec{r}_{ik}}{r_{ik}'+r_{ik}}\,,\nonumber \\[1ex]
\overline{\nabla}_{\vec{q}_j}^s V(\vec{q},\vec{q}') &=
  -\Delta_{r_{ji}}^s V(\vec{q},\vec{q}')
  \cdot
  \frac{\vec{r}_{ji}'+\vec{r}_{ji}}{r_{ji}'+r_{ji}}
  -
  \Delta_{r_{jk}}^s V(\vec{q},\vec{q}')
  \cdot
  \frac{\vec{r}_{jk}'+\vec{r}_{jk}}{r_{jk}'+r_{jk}} \,,\label{dgangles} \\[1ex]
\overline{\nabla}_{\vec{q}_k}^s V(\vec{q},\vec{q}') &=
  \phantom{-} \Delta_{r_{jk}}^s V(\vec{q},\vec{q}')
  \cdot
  \frac{\vec{r}_{jk}'+\vec{r}_{jk}}{r_{jk}'+r_{jk}}
  +
  \Delta_{r_{ik}}^s V(\vec{q},\vec{q}')
  \cdot
  \frac{\vec{r}_{ik}'+\vec{r}_{ik}}{r_{ik}'+r_{ik}} \,.\nonumber
\end{align}
These three expressions are indeed discrete gradients, which is detailed in the following theorem.
\begin{Thm} \label{thm:dgangle}
The three expressions are discrete gradients.
\end{Thm}
\begin{proof}
A calculation shows that the energy is preserved by the first gradient $\overline{\nabla}^\ell V$. We omit the arguments of the finite differences and obtain
\begin{align*}
 \overline{\nabla}^\ell &V(\vec{q},\vec{q}')^T (\vec{q}'-\vec{q}) \\
  &=
 \overline{\nabla}^\ell_{\vec{q}_i}V(\vec{q},\vec{q}')^T(\vec{q}_i'-\vec{q}_i) +
 \overline{\nabla}^\ell_{\vec{q}_j}V(\vec{q},\vec{q}')^T(\vec{q}_j'-\vec{q}_j) +
 \overline{\nabla}^\ell_{\vec{q}_k}V(\vec{q},\vec{q}')^T(\vec{q}_k'-\vec{q}_k) \\
  &=
  \Delta_{r_{ji}}^\ell V
  \cdot
  \left(\frac{\vec{r}_{ji}'+\vec{r}_{ji}}{r_{ji}'+r_{ji}}\right)^T
  (\vec{q}_i'-\vec{q}_i)
  -
  \Delta_{r_{ik}}^\ell V
  \cdot
  \left(\frac{\vec{r}_{ik}'+\vec{r}_{ik}}{r_{ik}'+r_{ik}}\right)^T
  (\vec{q}_i'-\vec{q}_i)
  \\
  &-
  \Delta_{r_{ji}}^\ell V
  \cdot
  \left(
  \frac{\vec{r}_{ji}'+\vec{r}_{ji}}{r_{ji}'+r_{ji}}
  \right)^T
  (\vec{q}_j'-\vec{q}_j)
  -
  \Delta_{r_{jk}}^\ell V
  \cdot
  \left(
  \frac{\vec{r}_{jk}'+\vec{r}_{jk}}{r_{jk}'+r_{jk}}
  \right)^T
  (\vec{q}_j'-\vec{q}_j) \\
  &+
  \Delta_{r_{jk}}^\ell V
  \cdot
  \left(
  \frac{\vec{r}_{jk}'+\vec{r}_{jk}}{r_{jk}'+r_{jk}}
  \right)^T
  (\vec{q}_k'-\vec{q}_k)
  +
  \Delta_{r_{ik}}^\ell V
  \cdot
  \left(
  \frac{\vec{r}_{ik}'+\vec{r}_{ik}}{r_{ik}'+r_{ik}}
  \right)^T
  (\vec{q}_k'-\vec{q}_k)
  \\
  &=
  \Delta_{r_{ji}}^\ell V
  \cdot
  \left(\frac{\vec{r}_{ji}'+\vec{r}_{ji}}{r_{ji}'+r_{ji}}\right)^T
  (\vec{r}_{ji}'-\vec{r}_{ji})
  +
  \Delta_{r_{jk}}^\ell V
  \cdot
  \left(
  \frac{\vec{r}_{jk}'+\vec{r}_{jk}}{r_{jk}'+r_{jk}}
  \right)^T
  (\vec{r}_{jk}'-\vec{r}_{jk})
  +
  \Delta_{r_{ik}}^\ell V
  \cdot
  \left(
  \frac{\vec{r}_{ik}'+\vec{r}_{ik}}{r_{ik}'+r_{ik}}
  \right)^T
  (\vec{r}_{ik}'-\vec{r}_{ik})
  \\
  &= \Delta_{r_{ji}}^\ell V \cdot (r_{ji}'-r_{ji})
   + \Delta_{r_{jk}}^\ell V \cdot (r_{jk}'-r_{jk})
   + \Delta_{r_{ik}}^\ell V \cdot (r_{ik}'-r_{ik}) \\
  &= V(r_{ji}',r_{jk},r_{ik})-V(r_{ji},r_{jk},r_{ik})
   + V(r_{ji}',r_{jk}',r_{ik})-V(r_{ji}',r_{jk},r_{ik}) 
   + V(r_{ji}',r_{jk}',r_{ik}')-V(r_{ji}',r_{jk}',r_{ik}) \\
  &= V(r_{ji}',r_{jk}',r_{ik}') - V(r_{ji},r_{jk},r_{ik})
   = V(\vec{q}')-V(\vec{q})\,.
\end{align*}
The discrete gradient is defined by continuous continuation. We have to show that the discrete gradient coincides with the exact gradient of $V$ in this case, that is, we have to show
\[
  \overline{\nabla}^\ell V(\vec{q},\vec{q})=\lim_{(\tilde{\vec{q}},\vec{q}') \rightarrow (\vec{q},\vec{q})} \overline{\nabla}^\ell V(\tilde{\vec{q}},\vec{q}') \stackrel{!}{=} \nabla V(\vec{q})\,.
\]
This is a simple calculation. We obtain
\begin{align*}
  \lim_{(\tilde{\vec{q}},\vec{q}')\rightarrow (\vec{q},\vec{q})} \overline{\nabla}^\ell_{\vec{q}_i}V(\tilde{\vec{q}},\vec{q}')
&=\phantom{-}
V_{r_{ji}}(r_{ji},r_{jk},r_{ik})\frac{\vec{r}_{ji}}{r_{ji}}
  -
V_{r_{ik}}(r_{ji},r_{jk},r_{ik})\frac{\vec{r}_{ik}}{r_{ik}}
=\nabla_{\vec{q}_i}V(\vec{q})\,,\\
  \lim_{(\tilde{\vec{q}},\vec{q}') \rightarrow (\vec{q},\vec{q})} \overline{\nabla}^\ell_{\vec{q}_j}V(\tilde{\vec{q}},\vec{q}')
&=
-V_{r_{ji}}(r_{ji},r_{jk},r_{ik})\frac{\vec{r}_{ji}}{r_{ji}}
  -
V_{r_{jk}}(r_{ji},r_{jk},r_{ik})\frac{\vec{r}_{jk}}{r_{jk}}
=\nabla_{\vec{q}_j}V(\vec{q})\,,\\
  \lim_{(\tilde{\vec{q}},\vec{q}') \rightarrow (\vec{q},\vec{q})} \overline{\nabla}^\ell_{\vec{q}_k}V(\tilde{\vec{q}},\vec{q}')
&=\phantom{-}
V_{r_{jk}}(r_{ji},r_{jk},r_{ik})\frac{\vec{r}_{jk}}{r_{jk}}
  +
V_{r_{ik}}(r_{ji},r_{jk},r_{ik})\frac{\vec{r}_{ik}}{r_{ik}}
=\nabla_{\vec{q}_k}V(\vec{q})\,.\\
\end{align*}
The proof of the discrete gradient $\overline{\nabla}^rV$ works analogously. From these two, finally the statement for the symmetric discrete gradient $\overline{\nabla}^sV$ follows.
\end{proof}
\noindent There are more discrete gradients. One might prescribe any pattern of primes to the three distances in the finite differences. For example the pattern prime, no prime, prime and then changing from prime to no prime and vice versa from left to right
\begin{align*}
\Delta_{r_{ji}}^f V(\vec{q},\vec{q}')
&\coloneqq
\frac{
   V(r_{ji}',r_{jk},r_{ik}')-V(r_{ji},r_{jk},r_{ik}')
  }{
   r_{ji}'-r_{ji}
  }\,,\\
\Delta_{r_{jk}}^f V(\vec{q},\vec{q}')
&\coloneqq
\frac{
   V(r_{ji}',r_{jk}',r_{ik}')-V(r_{ji}',r_{jk},r_{ik}')
  }{
   r_{jk}'-r_{jk}
  }\,,\\
\Delta_{r_{ik}}^f V(\vec{q},\vec{q}')
&\coloneqq
 \frac{
   V(r_{ji}',r_{jk},r_{ik}')-V(r_{ji}',r_{jk},r_{ik})
  }{
   r_{ik}'-r_{ik}
  }\,.
\end{align*}
All discrete gradients of this type can be symmetrised by symmetrising the finite differences.
\begin{Thm} \label{thm4}
 The method~\eqref{velocityDGstandardT} for a particle system of $N$ particles with Hamiltonian
\[
 H(\vec{q},\vec{p})=\frac{1}{2}\vec{p}^TM^{-1}\vec{p} + V(\vec{q})\,, \qquad V=\sum_{k=1}^M V^k\,,
\]	
where $V^k$ are angle potentials, is a first-order non-symmetric implicit method for $\overline{\nabla}^\ell V$ and $\overline{\nabla}^r V$ and a second-order symmetric implicit method for $\overline{\nabla}^s V$. All DG methods based on the discrete gradients above conserve the energy, the total linear momentum
	$\vec{P}=\sum_{i=1}^N \vec{p}_i$ and the total angular momentum $\vec{L}=\sum_{i=1}^N \vec{q}_i \times \vec{p}_i$.
\end{Thm}
\begin{proof}
The order of the methods follows from proposition~\ref{prop2}.
The preservation of the energy again follows from theorem~\ref{thm:pwdg}, lemma~\ref{lem:sumdg}, and lemma~\ref{energypressepHamil}. 	
From \eqref{dgangles}, or the corresponding relation for the other two discrete gradients, we have 
\[
  \sum_{i=1}^N \overline{\nabla}V^k_{\vec{q}_i} (\vec{q},\vec{q}')=0\,, \qquad \mbox{for} \qquad k=1,\ldots,M.	
\]
Hence, the preservation of the angular momentum follows literally as in equation~\eqref{thm2:eq1}.
For a fixed angle potential $V^k(\theta_{ijk})$ for particles $i$, $j$, $k$, for all three discrete gradients and for arbitrary $\vec{q}'$ and $\vec{q}$, we have that
\begin{align*}
&\sum_{i=1}^N (\vec{q}_i'+\vec{q}_i) \times \overline{\nabla}_{\vec{q}_i}V^k(\vec{q},\vec{q}')\\
   &=
   (\vec{q}_i'+\vec{q}_i)\times \overline{\nabla}_{\vec{q}_i} V^k(\vec{q},\vec{q}') +
   (\vec{q}_j'+\vec{q}_j)\times \overline{\nabla}_{\vec{q}_j} V^k(\vec{q},\vec{q}') + 
   (\vec{q}_k'+\vec{q}_k)\times \overline{\nabla}_{\vec{q}_k} V^k(\vec{q},\vec{q}') \\
  &= \left[ (\vec{q}_i'+\vec{q}_i) - (\vec{q}_j'+\vec{q}_j) \right] \times \Delta_{r_{ji}}V^k (\vec{q},\vec{q}')\cdot \frac{\vec{r}_{ji}'+\vec{r}_{ji}}{r_{ji}'+r_{ji}} 
   + \left[ (\vec{q}_k'+\vec{q}_k) - (\vec{q}_j'+\vec{q}_j) \right] \times \Delta_{r_{jk}}V^k (\vec{q},\vec{q}')\cdot \frac{\vec{r}_{jk}'+\vec{r}_{jk}}{r_{jk}'+r_{jk}} \\
  &+ \left[ (\vec{q}_k'+\vec{q}_k) - (\vec{q}_i'+\vec{q}_i) \right] \times \Delta_{r_{ik}}V^k (\vec{q},\vec{q}')\cdot \frac{\vec{r}_{ik}'+\vec{r}_{ik}}{r_{ik}'+r_{ik}} = \vec{0}\,.
\end{align*}
From here, the proof literally proceeds as in theorem~\ref{thm2}. 
\end{proof}
\noindent Theorem~\ref{thm2} and theorem~\ref{thm4} prove that the statement of theorem~\ref{thm4} remains true if pairwise potentials are included. This is summarised in our main theorem~\ref{thm:main}. The preservation of the total linear momentum and the total angular momentum is basically due to the following observation.
For the discrete forces
\[
  \overline{\vec{f}}_i=-\overline{\nabla}_{\vec{x}_i}V(\vec{q}^n,\vec{q}^{n+1}), \qquad 
  \overline{\vec{f}}_j=-\overline{\nabla}_{\vec{x}_j}V(\vec{q}^n,\vec{q}^{n+1}), \qquad 
  \overline{\vec{f}}_k=-\overline{\nabla}_{\vec{x}_k}V(\vec{q}^n,\vec{q}^{n+1}) 
\]
we have the relation analogue to the forces,
\[
\overline{\vec{f}}_k+\overline{\vec{f}}_i+\overline{\vec{f}}_j=0\,.    
\]
This is also satisfied by many other schemes, particularly the Verlet scheme. That is, these schemes also preserve the total linear momentum and the angular momentum. For this reason, we mainly consider the order and the energy preservation in our experiments. The energy preservation is what comes with the DG schemes.

\subsection{Experiment with bond angles}
\begin{figure}
\begin{center}
 \begin{minipage}{0.49\textwidth}	
\begin{tikzpicture}

\begin{axis}[
width=8cm, height=6.5cm, xmin=0.5*1e-6, xmax=2e-2, ymin=1e-10, ymax=40, xmode=log, ymode=log, legend pos=south east]
\addplot[mark=*, mark size=1pt, color=red!70!white, thick] coordinates {
( 1.000000000000000021e-02 , 1.718131008155034278e+00 )
( 1.000000000000000021e-03 , 1.724303209516057922e-02 )
( 1.000000000000000048e-04 , 1.723900479434825118e-04 )
( 1.000000000000000082e-05 , 1.723874723068426522e-06 )
( 9.999999999999999547e-07 , 1.719620985689127360e-08 )
};
\addplot[mark=diamond*, mark size = 2.5pt, thick, color=orange] coordinates {
( 1.000000000000000021e-02 , 2.757480375380190196e+00 )
( 1.000000000000000021e-03 , 3.276635217714132103e-02 )
( 1.000000000000000048e-04 , 3.278631181778646989e-04 )
( 1.000000000000000082e-05 , 3.278910994249723448e-06 )
( 9.999999999999999547e-07 , 3.298242856259716061e-08 )
};
\addplot[mark=+, mark size=2pt, color=cyan, thick] coordinates {
( 1.000000000000000021e-02 , 2.876092695910979824e+00 )
( 1.000000000000000021e-03 , 3.401634950103212679e-02 )
( 1.000000000000000048e-04 , 3.403612510743554416e-04 )
( 1.000000000000000082e-05 , 3.403653708286390223e-06 )
( 9.999999999999999547e-07 , 3.401490943883935718e-08 )
};
\addplot[color=black, thick, domain=0.000001:0.01] {100000*x*x};
\addplot[mark=*, color=blue!70!white, thick] coordinates {
( 1.000000000000000021e-02 , 4.083389775293396085e+00 )
( 1.000000000000000021e-03 , 6.945726505122541417e-01 )
( 1.000000000000000048e-04 , 6.589389592070012125e-02 )
( 1.000000000000000082e-05 , 6.566297022016135028e-03 )
( 9.999999999999999547e-07 , 6.564144395990736350e-04 )
};
\addplot[color=black, dashed, thick, domain=0.000001:0.01] {10000*x};

\addplot[color=black, thick, domain=0.000001:0.01] {100000*x*x};
\legend{Verlet, symmetric DG, midpoint, order2 line, simple DG, order1 line};

\end{axis}

\end{tikzpicture}
 \end{minipage}
 \begin{minipage}{0.49\textwidth}	
\input{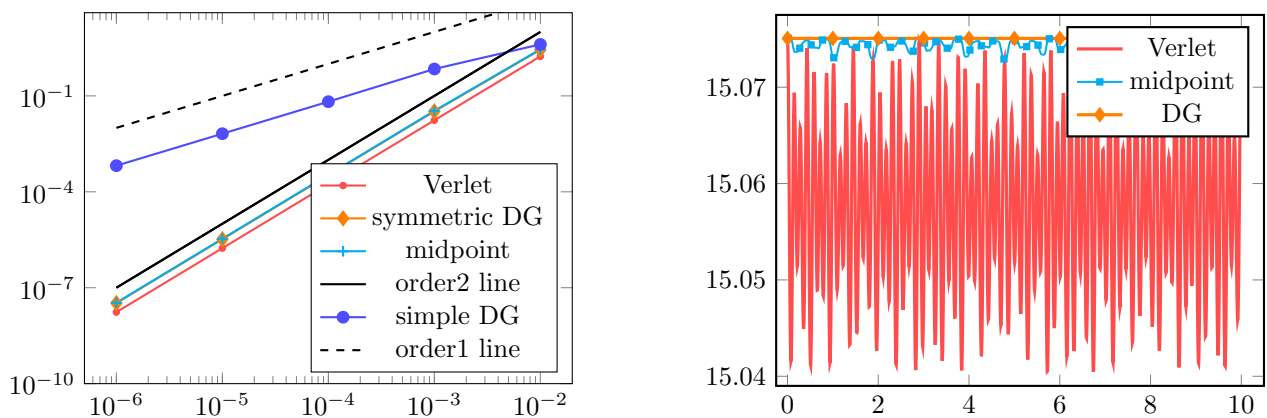}
 \end{minipage}
\end{center}
\caption{Results of the experiment with two water-like molecules: the error versus the time step is shown on the left-hand side for the discrete gradient (DG) methods, the midpoint rule, and the Verlet scheme. On the right-hand side, the energy is shown over the time span $[0,10]$ for step-size $\tau=0.005$ for all methods.} \label{exp:dgangle}
\end{figure}

We study the discrete gradients discussed so far by an experiment with molecules motivated by the TIP3P-C water model without long-range interactions. More exactly,
\[
  \begin{array}{rcc}
          \mbox{bond potential:} & r_0=0.957 \, \mbox{\AA} & k_b=450\,\frac{\mbox{\small kcal}}{\mbox{\small mol}} \\[2ex]
          \mbox{angle potential:} & \theta_0=104.52\degree & k_\theta = 55\,\frac{\mbox{\small kcal}}{\mbox{\small mol}} \\[2ex]
    \mbox{Lennard--Jones potential:} 
          & \epsilon_H=0.046\,\frac{\mbox{\small kcal}}{\mbox{\small mol}} & \sigma_H=0.4 \,\mbox{\AA} \\[1ex]
          & \epsilon_O=0.1521\,\frac{\mbox{\small kcal}}{\mbox{\small mol}} & \sigma_O=3.1506 \,\mbox{\AA} \\[1ex]
          & m_H=1.0080\,u & m_O=15.9994\,u 
  \end{array}
\]
We use the mixing rule
\[
  \sigma_{ij}=\frac{\sigma_{ii}+\sigma_{jj}}{2}, \qquad \epsilon_{ij}=\sqrt{\epsilon_{ii}\epsilon_{jj}} 
\]
whenever two different species $i$ and $j$ meet in a Lennard--Jones interaction. We only impose Lennard--Jones interactions on two particles that belong to different molecules. 
Here and in further experiments, we use the angle potential
\begin{equation} \label{anglepotcossq}
 V(\theta)=k_\theta(\cos(\theta)-\cos(\theta_0))^2\,.	
\end{equation}
This potential is available in LAMMPS, cf. \cite{Lammps22}, as {\tt angle\_style cosine/squared}.
Further, we remove the units by scaling. We use $\tilde{\sigma}=10^{-10}~\mbox{m}=\mbox{\AA}$, $\tilde{\epsilon}=1~\mbox{\small \mbox{kcal}}\slash\mbox{\small mol}$, $\tilde{m}=1~\mbox{u}$ and $\tilde{\alpha}=\tilde{\sigma}\sqrt{\tilde{m}\slash\tilde{\epsilon}}=4.88\cdot 10^{-14}$.
The dimensionless quantities then read
\[
 \begin{array}{lllll}
   m'=m\slash \tilde{m},~~~ & \vec{x}_i'=\vec{x}_i\slash \tilde{\sigma},~~~
   & \vec{r}_{ij}'=\vec{r}_{ij}\slash \tilde{\sigma},~~~
   & E'=E\slash \tilde{\epsilon},~~~ & V'=V\slash \tilde{\epsilon}, \\[1ex]
   \sigma'=\sigma\slash \tilde{\sigma},~~~ & \epsilon'=\epsilon\slash \tilde{\epsilon},~~~
         & t'=t\slash \tilde{\alpha}.~~~ & &
 \end{array}
\]
\begin{figure}
\begin{center}
\begin{minipage}{0.49\textwidth}
\renewcommand\lstlistingname{Code fragment}
\begin{lstlisting}[caption={Data}, escapechar=\#, label={dgangle}]
Positions

1 1  6.310593  9.000000 9.744371
2 1  7.054964  9.000000 9.000000
3 1  6.310593  9.000000 8.255629
4 2 11.689407  9.000000 9.744371
5 2 10.945036  9.000000 9.000000
6 2 11.689407  9.000000 8.255629

Velocities

1 0.000000 -0.050000 0.000000
2 0.000000 -0.050000 0.000000
3 0.000000 -0.050000 0.000000
4 0.000000  0.050000 0.000000
5 0.000000  0.050000 0.000000
6 0.000000  0.050000 0.000000
\end{lstlisting}
\end{minipage}
\hspace{-2cm}
\begin{minipage}{0.49\textwidth}
\begin{center}
\hspace*{1cm} \includegraphics[height=3cm]{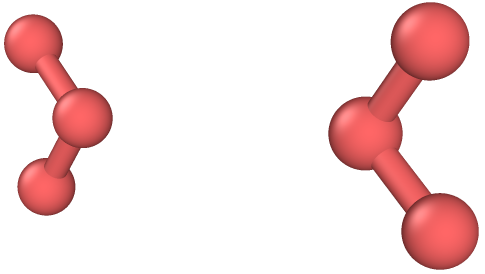}\\[4ex]
\end{center}
\end{minipage}
\end{center}
\caption{Initial conditions for the experiment with two water-like molecules: the positions and velocities are given on the left-hand side. In the positions section, the first row numbers the particles. The second row assigns molecule numbers. If the number is the same, a bond is added between the particles. The third to fifth row are the coordinates of the positions. In the velocities section, the particle number is followed by velocities. On the right-hand side, the initial configuration is given as a plot by the open visualization tool (OVITO), cf.\cite{Stukowski10}} \label{dgangleinitial}
\end{figure}
In figure~\ref{dgangleinitial}, the (dimensionless) initial positions and velocities of the particles are given. The first line in the positions section numbers the atoms. The second indicates the molecule number. As long as this number is the same, the atoms are connected by a bond. The velocities of the particles are given in the next section. We first compute the trajectory up to time $10$. We show the error versus the step size at time $10$ on the left-hand side of figure~\ref{exp:dgangle}. All methods show their expected order.
With an unsymmetric discrete gradient, we observe first order. When the discrete gradient is symmetrised, second order for the DG scheme is observed.
On the right-hand side of figure~\ref{exp:dgangle}, the calculated energies for time step $\tau=0.005$ are shown versus the simulation time. All three discrete gradient schemes preserve the energy nicely.
The implicit equation in \eqref{velocityDGmethod} is solved with the simplified Newton method, \eqref{SimpleNewtonmeth}, where the approximate Jacobian $J_F$, cf. \eqref{simpleJac}, includes the full Hessian with respect to the Lennard--Jones potentials and the bond potentials but omits the part with respect to the bond angle potentials. 
\subsection{Discrete gradients for dihedral angles}
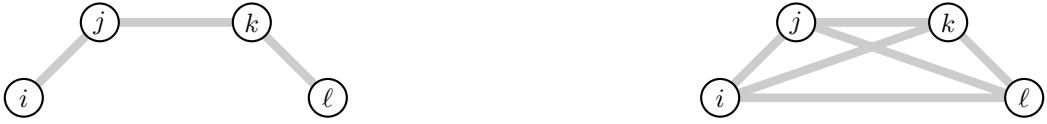
\begin{figure}[h]
\begin{center}
\begin{minipage}{0.49\textwidth}
\begin{center}
\begin{tikzpicture}[scale=1.0]
         \filldraw[white] (1,1) circle (0.25cm);
         \filldraw[white] (2,2) circle (0.25cm);
         \filldraw[white] (4,2) circle (0.25cm);
         \filldraw[white] (5,1) circle (0.25cm);
         \pgfmathsetmacro\xic{1+0.25*cos(atan(1))};
         \pgfmathsetmacro\xis{1+0.25*sin(atan(1))};
         \pgfmathsetmacro\xjcl{2-0.25*cos(atan(1))};
         \pgfmathsetmacro\xjsl{2-0.25*sin(atan(1))};
         \pgfmathsetmacro\xkcr{4+0.25*cos(atan(1))};
         \pgfmathsetmacro\xksr{2-0.25*sin(atan(1))};
         \pgfmathsetmacro\xlc{5-0.25*cos(atan(1))};
         \pgfmathsetmacro\xls{1+0.25*sin(atan(1))};
         \pgfmathsetmacro\xicoff{0.05/sqrt(2)};
         \pgfmathsetmacro\xisoff{0.05/sqrt(2)};
         \filldraw[white!80!black] (\xic+\xicoff,\xis-\xisoff) -- (\xjcl+\xicoff,\xjsl-\xisoff) -- (\xjcl-\xicoff,\xjsl+\xisoff) -- (\xic-\xicoff,\xis+\xisoff);
         \filldraw[white!80!black] (\xkcr+\xicoff,\xksr+\xisoff) -- (\xlc+\xicoff,\xls+\xisoff) -- (\xlc-\xicoff,\xls-\xisoff) -- (\xkcr-\xicoff,\xksr-\xisoff);
         \filldraw[white!80!black] (2+0.25,2+0.05) -- (4-0.25,2+0.05) -- (4-0.25,2-0.05) -- (2+0.25,2-0.05);
         \draw[thick] (1,1) circle (0.25cm);
         \node at (1,1) {$i$};
         \draw[thick] (2,2) circle (0.25cm);
         \node at (2,2) {$j$};
         \draw[thick] (4,2) circle (0.25cm);
         \node at (4,2) {$k$};
         \draw[thick] (5,1) circle (0.25cm);
         \node at (5,1) {$\ell$};
\end{tikzpicture}
\end{center}
\end{minipage}
\begin{minipage}{0.49\textwidth}
\begin{center}
\begin{tikzpicture}[scale=1.0]
         \filldraw[white] (1,1) circle (0.25cm);
         \filldraw[white] (2,2) circle (0.25cm);
         \filldraw[white] (4,2) circle (0.25cm);
         \filldraw[white] (5,1) circle (0.25cm);
         \pgfmathsetmacro\xic{1+0.25*cos(atan(1))};
         \pgfmathsetmacro\xis{1+0.25*sin(atan(1))};
         \pgfmathsetmacro\xjcl{2-0.25*cos(atan(1))};
         \pgfmathsetmacro\xjsl{2-0.25*sin(atan(1))};
         \pgfmathsetmacro\xkcr{4+0.25*cos(atan(1))};
         \pgfmathsetmacro\xksr{2-0.25*sin(atan(1))};
         \pgfmathsetmacro\xlc{5-0.25*cos(atan(1))};
         \pgfmathsetmacro\xls{1+0.25*sin(atan(1))};
         \pgfmathsetmacro\xicoff{0.05/sqrt(2)};
         \pgfmathsetmacro\xisoff{0.05/sqrt(2)};
         \filldraw[white!80!black] (\xic+\xicoff,\xis-\xisoff) -- (\xjcl+\xicoff,\xjsl-\xisoff) -- (\xjcl-\xicoff,\xjsl+\xisoff) -- (\xic-\xicoff,\xis+\xisoff);
         \filldraw[white!80!black] (\xkcr+\xicoff,\xksr+\xisoff) -- (\xlc+\xicoff,\xls+\xisoff) -- (\xlc-\xicoff,\xls-\xisoff) -- (\xkcr-\xicoff,\xksr-\xisoff);
	 \filldraw[white!80!black] (2+0.25,2+0.05) -- (4-0.25,2+0.05) -- (4-0.25,2-0.05) -- (2+0.25,2-0.05);
	 \pgfmathsetmacro\xikc{1+0.25*cos(atan(0.25))};
	 \pgfmathsetmacro\xiks{1+0.25*sin(atan(0.25))};
	 \pgfmathsetmacro\xkic{4-0.25*cos(atan(0.25))};
	 \pgfmathsetmacro\xkis{2-0.25*sin(atan(0.25))};
         \pgfmathsetmacro\xicoff{-0.05*sin(atan(0.25))};
         \pgfmathsetmacro\xisoff{0.05*cos(atan(0.25))};
	 \filldraw[white!80!black] (\xikc+\xicoff,\xiks+\xisoff) -- (\xkic+\xicoff,\xkis+\xisoff) -- (\xkic-\xicoff,\xkis-\xisoff) -- (\xikc-\xicoff,\xiks-\xisoff);
	 \pgfmathsetmacro\xjlc{2+0.25*cos(atan(0.25))};
	 \pgfmathsetmacro\xjls{2-0.25*sin(atan(0.25))};
	 \pgfmathsetmacro\xljc{5-0.25*cos(atan(0.25))};
	 \pgfmathsetmacro\xljs{1+0.25*sin(atan(0.25))};
         \pgfmathsetmacro\xicoff{0.05*sin(atan(0.25))};
         \pgfmathsetmacro\xisoff{0.05*cos(atan(0.25))};
	 \filldraw[white!80!black] (\xjlc+\xicoff,\xjls+\xisoff) -- (\xljc+\xicoff,\xljs+\xisoff) -- (\xljc-\xicoff,\xljs-\xisoff) -- (\xjlc-\xicoff,\xjls-\xisoff);
	 \filldraw[white!80!black] (1+0.25,1+0.05) -- (5-0.25,1+0.05) -- (5-0.25,1-0.05) -- (1+0.25,1-0.05);
         \draw[thick] (1,1) circle (0.25cm);
         \node at (1,1) {$i$};
         \draw[thick] (2,2) circle (0.25cm);
         \node at (2,2) {$j$};
         \draw[thick] (4,2) circle (0.25cm);
         \node at (4,2) {$k$};
         \draw[thick] (5,1) circle (0.25cm);
         \node at (5,1) {$\ell$};
\end{tikzpicture}
\end{center}
\end{minipage}
\end{center}
\caption{Sketch of united-atom butane: on the left-hand side the molecule with the particles $i$, $j$, $k$, $\ell$ is shown including the bonds between atoms. On the right-hand side, the distances between the particles $i$, $j$, $k$, $\ell$ used by the discrete gradient are given.} \label{dihedralplot}
\end{figure}
\noindent The dihedral angle
\begin{align} \label{dihedralIUPAC}
        \phi_{ijk\ell} &= \sign( \langle \vec{r}_{ij}, \vec{r}_{jk}\times \vec{r}_{k\ell}\rangle) \arccos \varphi_{ijk\ell}, \qquad \varphi_{ijk\ell}=
   \frac{
     \langle
       \vec{r}_{ij} \times \vec{r}_{jk}, \vec{r}_{jk} \times \vec{r}_{k\ell}
     \rangle
    }{
     \| \vec{r}_{ij} \times \vec{r}_{jk} \| \|\vec{r}_{jk} \times \vec{r}_{k\ell}\|
    }\,,
\end{align}
denotes the angle between the planes spanned by the atoms $i$, $j$, $k$ and $j$, $k$, $\ell$, respectively. The sign of the dihedral angle $\phi_{ijk\ell}$, i.e. $\sign( \langle \vec{r}_{ij}, \vec{r}_{jk}\times \vec{r}_{k\ell} \rangle)$, designates on which side of the plane through $j$, $k$, $\ell$ the particle $i$ lies. We consider potentials $V(\phi)$ that depend smoothly on the dihedral angle $\phi$. The computation of the forces based on the dihedral angle is based on the expression as given above in \eqref{dihedralIUPAC} and leads to the most used formulas as described in \cite{Bekker96}. To our best knowledge, the discrete gradient that we propose in the following is new. It is based on the same idea as the discrete gradients before, namely, the idea to express the angle with distances. The dihedral angle can also be expressed in terms of the distances between all atoms as indicated on the right-hand side of figure~\ref{dihedralplot}. The additional lines on the right-hand side of figure~\ref{dihedralplot} are the additional distances that will be used. The representation of the angle in distances is given in the following lemma.
\begin{lemma} \label{lem:dihedraldist}
We have the following representation in terms of distances
\begin{align*}
  \varphi_{ijk\ell} &=
  \frac{
    \langle
        \vec{m}, \vec{n}
    \rangle
  }{
          \| \vec{m} \| \|\vec{n}\|
  }
=
\frac{
 (r_{k\ell}^2+r_{jk}^2-r_{j\ell}^2)(r_{ij}^2+r_{jk}^2-r_{ik}^2)
 -2r_{jk}^2(r_{jk}^2+r_{i\ell}^2-r_{j\ell}^2-r_{ik}^2)
}{
\sqrt{4r_{jk}^2r_{ij}^2-(r_{ij}^2+r_{jk}^2-r_{ik}^2)^2}
\sqrt{4r_{jk}^2r_{k\ell}^2-(r_{jk}^2+r_{k\ell}^2-r_{j\ell}^2)^2}
}
\end{align*}
where $\vec{m} = \vec{r}_{ij} \times \vec{r}_{jk}$ and $\vec{n} = \vec{r}_{jk} \times \vec{r}_{k\ell}$.
\end{lemma}
\begin{proof}
The proof is a tedious calculation along the following steps
\begin{align*}
 \frac{
    \langle
        \vec{r}_{ij} \times \vec{r}_{jk}, \vec{r}_{jk} \times \vec{r}_{k\ell}
    \rangle
  }{
   \| \vec{r}_{ij} \times \vec{r}_{jk} \| \|\vec{r}_{jk} \times \vec{r}_{k\ell}\|
  }
&=
(-1)\cdot
  \frac{
    \Big\langle
      r_{jk}^2\vec{r}_{ij} - \big\langle \vec{r}_{ij}, \vec{r}_{jk} \big\rangle \vec{r}_{jk},
      r_{jk}^2\vec{r}_{k\ell} - \big\langle \vec{r}_{k\ell}, \vec{r}_{jk}\big\rangle \vec{r}_{jk}
    \Big\rangle
  }{
   \|
     r_{jk}^2\vec{r}_{ij} - \big\langle \vec{r}_{ij}, \vec{r}_{jk}\big\rangle \vec{r}_{jk}
   \|
   \|
    r_{jk}^2\vec{r}_{k\ell} - \big\langle \vec{r}_{k\ell}, \vec{r}_{jk}  \big\rangle \vec{r}_{jk}
   \|} \\[1.5ex]
&=
\frac{
 \langle \vec{r}_{k\ell}, \vec{r}_{jk}\rangle
 \langle \vec{r}_{ij}, \vec{r}_{jk} \rangle
 - r_{jk}^2 \langle \vec{r}_{ij}, \vec{r}_{k\ell} \rangle
}{
\sqrt{r_{jk}^2r_{ij}^2 - \langle \vec{r}_{ij}, \vec{r}_{jk}\rangle^2}
\sqrt{r_{jk}^2r_{k\ell}^2 - \langle \vec{r}_{jk}, \vec{r}_{k\ell}\rangle^2}
}\\[1.5ex]
&=
\frac{
 (r_{k\ell}^2+r_{jk}^2-r_{j\ell}^2)(r_{ij}^2+r_{jk}^2-r_{ik}^2)
 -2r_{jk}^2(r_{jk}^2+r_{i\ell}^2-r_{j\ell}^2-r_{ik}^2)
}{
\sqrt{4r_{jk}^2r_{ij}^2-(r_{ij}^2+r_{jk}^2-r_{ik}^2)^2}
\sqrt{4r_{jk}^2r_{k\ell}^2-(r_{jk}^2+r_{k\ell}^2-r_{j\ell}^2)^2}
}
\end{align*}
\end{proof}
\noindent So the dihedral potential depends on six distances $V(r_{ij},r_{jk},r_{k\ell},r_{ik},r_{j\ell},r_{i\ell})$.
The representation in lemma~\ref{lem:dihedraldist} suggests an alternative way to represent any potential dependent on the dihedral angle in terms of distances. But we will consider potentials that are dependent on the cosine of the dihedral angle, because we will use such potentials in our experiments later. Let
\[
        U_t(\phi)=k_\phi \sum_{n=0}^m a_n \cos^n \phi=\widetilde{U}_t(\cos \phi)\,, \qquad \widetilde{U}_t(\varphi)=k_\phi\sum_{n=0}^m a_n \varphi^n\,,
\]
be the potential for the torsion angles.
Then, we have
\begin{align*}
 \frac{\partial}{\partial \vec{q}_u} U_t(\phi_{ijk\ell})
  &=
        \tilde{U}_t'(\cos \phi_{ijk\ell})\cdot (-1)\cdot \sin(\phi_{ijk\ell}) \sign(\langle \vec{r}_{ij},\vec{n}\rangle) \frac{-1}{\sin \varphi_{ijk\ell}} \cdot \frac{\partial}{\partial{\vec{q}_u}} \varphi_{ijk\ell} \\
  &=
        \tilde{U}_t'(\varphi_{ijk\ell}) \cdot \frac{\partial}{\partial{\vec{q}_u}} \varphi_{ijk\ell}\,, \qquad u=i,j,k,\ell\,.
\end{align*}
The potential reduces to
$
 \tilde{U}_t(\varphi_{ijk\ell})
$.
Hence, for the discrete gradient, we use the representation
\[
 V(r_{ij},r_{jk},r_{k\ell},r_{ik},r_{j\ell},r_{i\ell})=\tilde{U}_t(\varphi_{ijk\ell}),
\]
with $\varphi_{ijk\ell}$ expressed in the distances as in lemma~\ref{lem:dihedraldist}.

\begin{figure}[h]
\begin{center}
\begin{minipage}{0.49\textwidth}
\begin{center}
\begin{tikzpicture}[scale=1.0]
         \filldraw[white] (1,1) circle (0.25cm);
         \filldraw[white] (2,2) circle (0.25cm);
         \filldraw[white] (4,2) circle (0.25cm);
         \filldraw[white] (5,1) circle (0.25cm);
         \pgfmathsetmacro\xk{2+sqrt(2)};
         \pgfmathsetmacro\xic{1+0.25*cos(atan(1))};
         \pgfmathsetmacro\xis{1+0.25*sin(atan(1))};
         \pgfmathsetmacro\xjcl{2-0.25*cos(atan(1))};
         \pgfmathsetmacro\xjsl{2-0.25*sin(atan(1))};
         \pgfmathsetmacro\xjlcr{2-0.25*cos(atan(1))};
         \pgfmathsetmacro\xjlsr{2+0.25*sin(atan(1))};
         \pgfmathsetmacro\xlc{1+0.25*cos(atan(1))};
         \pgfmathsetmacro\xls{3-0.25*sin(atan(1))};
         \pgfmathsetmacro\xicoff{0.05/sqrt(2)};
         \pgfmathsetmacro\xisoff{0.05/sqrt(2)};
         \filldraw[white!80!black] (\xic+\xicoff,\xis-\xisoff) -- (\xjcl+\xicoff,\xjsl-\xisoff) -- (\xjcl-\xicoff,\xjsl+\xisoff) -- (\xic-\xicoff,\xis+\xisoff);
         \filldraw[white!80!black] (\xjlcr+\xicoff,\xjlsr+\xisoff) -- (\xlc+\xicoff,\xls+\xisoff) -- (\xlc-\xicoff,\xls-\xisoff) -- (\xjlcr-\xicoff,\xjlsr-\xisoff);
         \filldraw[white!80!black] (2+0.25,2+0.05) -- (\xk-0.25,2+0.05) -- (\xk-0.25,2-0.05) -- (2+0.25,2-0.05);
         \draw[thick] (1,1) circle (0.25cm);
         \node at (1,1) {$i$};
         \draw[thick] (2,2) circle (0.25cm);
         \node at (2,2) {$k$};
         \draw[thick] (\xk,2) circle (0.25cm);
	 \node at (\xk,2) {$j$};
         \draw[thick] (1,3) circle (0.25cm);
         \node at (1,3) {$\ell$};
\end{tikzpicture}
\end{center}
\end{minipage}
\hspace*{-2.5cm}
\begin{minipage}{0.49\textwidth}
\[
	\omega_{ijk\ell} = \arccos(\varphi_{ijk\ell}), \qquad
\varphi_{ijk\ell}=
   \frac{
     \langle
       \vec{r}_{ij} \times \vec{r}_{jk}, \vec{r}_{jk} \times \vec{r}_{k\ell}
     \rangle
    }{
     \| \vec{r}_{ij} \times \vec{r}_{jk} \| \|\vec{r}_{jk} \times \vec{r}_{k\ell}\|
    }
\]
\end{minipage}
\end{center}
	\caption{Improper dihedral angle: on the left-hand side, the typical situation is shown. On the right-hand side, the formula for the improper dihedral angle $\omega_{ijk\ell}$ is given. $\varphi_{ijk\ell}$ is exactly the same term as for the standard dihedral angle, cf. \eqref{dihedralIUPAC}.}\label{improperplot}
\end{figure}
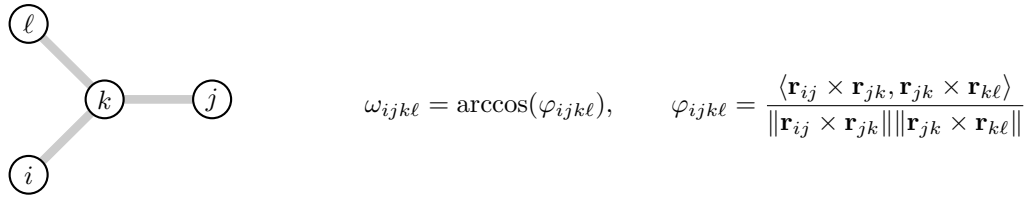

We first discuss non-symmetric discrete gradients. In order to write the discrete gradient down in a concise way, we use again finite differences with respect to the distances of the particles.
\[
\begin{array}{lcr}
\Delta_{r_{ij}}^\ell V(\vec{q},\vec{q}')
	&\coloneqq&
{\displaystyle
  \frac{
         V(r_{ij}',r_{jk},r_{k\ell},r_{ik},r_{j\ell},r_{i\ell})
        -V(r_{ij},r_{jk},r_{k\ell},r_{ik},r_{j\ell},r_{i\ell})
  }{
   r_{ij}'-r_{ij}
  }
}\\
&\vdots& \\
\Delta_{r_{i\ell}}^\ell V(\vec{q},\vec{q}')
&\coloneqq&
{\displaystyle
\frac{
        V(r_{ij}',r_{jk}',r_{k\ell}',r_{ik}',r_{j\ell}',r_{i\ell}')
        -V(r_{ij}',r_{jk}',r_{k\ell}',r_{ik}',r_{j\ell}',r_{i\ell})
  }{	  
   r_{i\ell}'-r_{i\ell}
  }
}	
\end{array}	
\]
This can also be seen as the Itoh--Abe (coordinate increment) discrete gradient for $V$. With these finite differences, our first discrete gradient reads as follows
\begin{equation} \label{dgdihedral}
\begin{array}{rcl}
  \overline{\nabla}_{\vec{q}_i}^\ell V(\vec{q},\vec{q}') &=&
{\displaystyle	
  -\Delta_{r_{ij}}^\ell V
  \cdot
  \frac{\vec{r}_{ij}'+\vec{r}_{ij}}{r_{ij}'+r_{ij}}
  -
  \Delta_{r_{ik}}^\ell V
  \cdot
  \frac{\vec{r}_{ik}'+\vec{r}_{ik}}{r_{ik}'+r_{ik}}
  -
  \Delta_{r_{i\ell}}^\ell V
  \cdot
  \frac{\vec{r}_{i\ell}'+\vec{r}_{i\ell}}{r_{i\ell}'+r_{i\ell}} 
}\,,\\[1.5ex]
  \overline{\nabla}_{\vec{q}_j}^\ell V(\vec{q},\vec{q}') &=&
{\displaystyle
  \phantom{-} \Delta_{r_{ij}}^\ell V
  \cdot
  \frac{\vec{r}_{ij}'+\vec{r}_{ij}}{r_{ij}'+r_{ij}}
  -
  \Delta_{r_{jk}}^\ell V
  \cdot
  \frac{\vec{r}_{jk}'+\vec{r}_{jk}}{r_{jk}'+r_{jk}}
  -
  \Delta_{r_{j\ell}}^\ell V
  \cdot
  \frac{\vec{r}_{j\ell}'+\vec{r}_{j\ell}}{r_{j\ell}'+r_{j\ell}} 
} \,,\\[1.5ex]
  \overline{\nabla}_{\vec{q}_k}^\ell V(\vec{q},\vec{q}') &=&
{\displaystyle
  \phantom{-}\Delta_{r_{ik}}^\ell V
  \cdot
  \frac{\vec{r}_{ik}'+\vec{r}_{ik}}{r_{ik}'+r_{ik}}
  +
  \Delta_{r_{jk}}^\ell V
  \cdot
  \frac{\vec{r}_{jk}'+\vec{r}_{jk}}{r_{jk}'+r_{jk}}
  -
  \Delta_{r_{k\ell}}^\ell V
  \cdot
  \frac{\vec{r}_{k\ell}'+\vec{r}_{k\ell}}{r_{k\ell}'+r_{k\ell}} 
} \,,\\[1.5ex]
  \overline{\nabla}_{\vec{q}_\ell}^\ell V(\vec{q},\vec{q}') &=&
{\displaystyle
  \phantom{-}\Delta_{r_{i\ell}}^\ell V
  \cdot
  \frac{\vec{r}_{i\ell}'+\vec{r}_{i\ell}}{r_{i\ell}'+r_{i\ell}}
  +
  \Delta_{r_{j\ell}}^\ell V
  \cdot
  \frac{\vec{r}_{j\ell}'+\vec{r}_{j\ell}}{r_{j\ell}'+r_{j\ell}}
  +
  \Delta_{r_{k\ell}}^\ell V
  \cdot
  \frac{\vec{r}_{k\ell}'+\vec{r}_{k\ell}}{r_{k\ell}'+r_{k\ell}}\,.
}
\end{array}
\end{equation}
This discrete gradient is not symmetric. The coordinate decrement version
\[
\begin{array}{rcl}
\Delta_{r_{ij}}^r V(\vec{q},\vec{q}')
&\coloneqq&
{\displaystyle
\frac{
        V(r_{ij},r_{jk}',r_{k\ell}',r_{ik}',r_{j\ell}',r_{i\ell}')
        -V(r_{ij}',r_{jk}',r_{k\ell}',r_{ik}',r_{j\ell}',r_{i\ell}')
  }{
        r_{ij}-r_{ij}'
  }
}\\
&\vdots& \\
\Delta_{r_{i\ell}}^r V(\vec{q},\vec{q}')
&\coloneqq&
{\displaystyle
\frac{
        V(r_{ij},r_{jk},r_{k\ell},r_{ik},r_{j\ell},r_{i\ell})
        -V(r_{ij},r_{jk},r_{k\ell},r_{ik},r_{j\ell},r_{i\ell}')
  }{
        r_{i\ell}-r_{i\ell}'
  }
}
\end{array}
\]
leads to another discrete gradient that is also not symmetric.
Since symmetric discrete gradients lead to symmetric time-integration methods, which are of second order, one might prefer symmetric discrete gradients. In the same way as the Itoh--Abe discrete gradient, we can symmetrise our gradients here. This leads to the discrete gradient with the (symmetric) finite differences
\begin{align*}
   \Delta_{r}^s V(\vec{q},\vec{q}')
   \coloneqq\frac{1}{2}
    \left(
        \Delta_{r}^\ell V(\vec{q},\vec{q}')+\Delta_{r}^r V(\vec{q},\vec{q}')
        \right), \quad r \in \{r_{ij}, r_{jk}, r_{kl}, r_{ik}, r_{j\ell}, r_{i\ell} \}
\end{align*}
These three expressions are indeed discrete gradients, which is detailed in the following theorem. 
The proof of this theorem is analogous to the proof of theorem~\ref{thm:dgangle}.
\begin{Thm}
The three expressions are discrete gradients.
\end{Thm}
\noindent Actually, there are many more discrete gradients. One might prescribe any pattern of primes to the six distances. Changing the primed distances to not-primed distances and vice versa from left to right when the finite difference 'passed' this distance, we obtain a new unsymmetric discrete gradient. These unsymmetric discrete gradients can all be symmetrised.
\begin{Thm} \label{thm6}
 The method~\eqref{velocityDGstandardT} for a particle system of $N$ particles with Hamiltonian
\[
 H(\vec{q},\vec{p})=\frac{1}{2}\vec{p}^TM^{-1}\vec{p} + V(\vec{q})\,, \qquad V=\sum_{k=1}^M V^k\,,
\]	
where $V^k$ are dihedral potentials, is a first-order non-symmetric implicit method for $\overline{\nabla}^\ell V$ and $\overline{\nabla}^r V$ and a second-order symmetric implicit method for $\overline{\nabla}^s V$. All three methods conserve the energy, the total linear momentum
$\vec{P}=\sum_{i=1}^N \vec{p}_i$ and the total angular momentum $\vec{L}=\sum_{i=1}^N \vec{q}_i \times \vec{p}_i$.
\end{Thm}
\begin{proof}
The order of the methods follows from proposition~\ref{prop2}.
The preservation of the energy again follows from theorem~\ref{thm:pwdg}, lemma~\ref{lem:sumdg}, and lemma~\ref{energypressepHamil}. 	
From \eqref{dgdihedral}, or the corresponding relation for the other two discrete gradients, we have 
\[
  \sum_{i=1}^N \overline{\nabla}V^k_{\vec{q}_i} (\vec{q},\vec{q}')=0\,, \qquad \mbox{for} \qquad k=1,\ldots,M.	
\]
Hence, the preservation of the angular momentum follows literally as in equation~\eqref{thm2:eq1}.
For a fixed dihedral potential $V^k(\varphi_{ijk\ell})$ for particles $i$, $j$, $k$, $\ell$, for all three discrete gradients and for arbitrary $\vec{q}'$ and $\vec{q}$, we have that
\begin{align*}
  \sum_{i=1}^N &(\vec{q}_i'+\vec{q}_i) \times \overline{\nabla}_{\vec{q}_i}V^k(\vec{q},\vec{q}') \\
  &=
   (\vec{q}_i'+\vec{q}_i)\times \overline{\nabla}_{\vec{q}_i} V^k(\vec{q},\vec{q}') +
   (\vec{q}_j'+\vec{q}_j)\times \overline{\nabla}_{\vec{q}_j} V^k(\vec{q},\vec{q}') + 
   (\vec{q}_k'+\vec{q}_k)\times \overline{\nabla}_{\vec{q}_k} V^k(\vec{q},\vec{q}') +
   (\vec{q}_\ell'+\vec{q}_\ell)\times \overline{\nabla}_{\vec{q}_\ell} V^k(\vec{q},\vec{q}')
	\\
  &= \left[ (\vec{q}_j'+\vec{q}_j) - (\vec{q}_i'+\vec{q}_i) \right] \times \Delta_{r_{ij}}V^k (\vec{q},\vec{q}')\cdot \frac{\vec{r}_{ij}'+\vec{r}_{ij}}{r_{ij}'+r_{ij}} 
   + \left[ (\vec{q}_k'+\vec{q}_k) - (\vec{q}_i'+\vec{q}_i) \right] \times \Delta_{r_{ik}}V^k (\vec{q},\vec{q}')\cdot \frac{\vec{r}_{ik}'+\vec{r}_{ik}}{r_{ik}'+r_{ik}} \\
  &+ \left[ (\vec{q}_\ell'+\vec{q}_\ell) - (\vec{q}_i'+\vec{q}_i) \right] \times \Delta_{r_{i\ell}}V^k (\vec{q},\vec{q}')\cdot \frac{\vec{r}_{i\ell}'+\vec{r}_{i\ell}}{r_{i\ell}'+r_{i\ell}} 
	+ \left[ (\vec{q}_k'+\vec{q}_k) - (\vec{q}_j'+\vec{q}_j) \right] \times \Delta_{r_{jk}}V^k (\vec{q},\vec{q}')\cdot \frac{\vec{r}_{jk}'+\vec{r}_{jk}}{r_{jk}'+r_{jk}}\\
  &+ \left[ (\vec{q}_\ell'+\vec{q}_\ell) - (\vec{q}_j'+\vec{q}_j) \right] \times \Delta_{r_{j\ell}}V^k (\vec{q},\vec{q}')\cdot \frac{\vec{r}_{j\ell}'+\vec{r}_{j\ell}}{r_{j\ell}'+r_{j\ell}} 
	+ \left[ (\vec{q}_\ell'+\vec{q}_\ell) - (\vec{q}_k'+\vec{q}_k) \right] \times \Delta_{r_{k\ell}}V^k (\vec{q},\vec{q}')\cdot \frac{\vec{r}_{k\ell}'+\vec{r}_{k\ell}}{r_{k\ell}'+r_{k\ell}} = \vec{0}\,.\\
\end{align*}
From here, the proof literally proceeds as in theorem~\ref{thm2}. 
\end{proof}
\noindent Besides the conventional torsion potential, where the atoms $i$, $j$, $k$, $\ell$ are consecutively connected, there is also the possibility that three atoms are connected to one as in figure~\ref{improperplot} on the left-hand side. Here, the so-called improper dihedral angle is defined as the angle between the planes spanned by the atoms $i$,$j$,$k$ and $j$,$k$,$\ell$, respectively. As shown on the right-hand side of figure~\ref{improperplot}, this angle can be expressed in exactly the same way as before for the standard dihedral angle \eqref{dihedralIUPAC}. Therefore, the distance definition of the dihedral angle in lemma~\ref{lem:dihedraldist} also works for the improper dihedral angle and the corresponding discrete gradients are constructed analogously.  
\noindent As a corollary from theorem~\ref{thm2}, theorem~\ref{thm4}, and theorem~\ref{thm6}, we obtain our main theorem.
\begin{Thm} \label{thm:main}
 The method~\eqref{velocityDGstandardT} for a particle system of $N$ particles with Hamiltonian
\[
 H(\vec{q},\vec{p})=\frac{1}{2}\vec{p}^TM^{-1}\vec{p} + V(\vec{q})\,, \qquad V=\sum_{k=1}^M V^k\,,
\]	
where $V^k$ are pairwise, angle, and dihedral potentials, is a first-order non-symmetric implicit method if at least one non-symmetric discrete gradient is used, and a second-order symmetric implicit method if all discrete gradients are symmetric. All methods preserve the energy, the total linear momentum
$\vec{P}=\sum_{i=1}^N \vec{p}_i$ and the total angular momentum $\vec{L}=\sum_{i=1}^N \vec{q}_i \times \vec{p}_i$.
\end{Thm}
\noindent Theorem~\ref{thm:main} shows that all standard short-range forces in a classical molecular dynamics simulation can be modelled with discrete gradients. Since even large protein simulations are often based on these standard short-range potentials and since the standard potentials do not change whether more atoms are connected to atoms defining one of the standard interactions or not, this theorem is quite general with respect to its applicability. 

Additional insight in the computation of the gradient of a dihedral angle potential might be gained from the representation of the dihedral angle with six distances. The computation of this gradient is not straightforward. It is usually based on formula~\eqref{dihedralIUPAC}, also called the cross-product definition of the dihedral angle. The scalar product of the cross-products can be represented as a scalar product without cross products, the so-called scalar product definition for the torsion angle. The cross-product definition leads to the most used formulas for the computation of the forces induced by the dihedral potential as given in \cite{Bekker96}. 
The representation of the dihedral angle in terms of distances leads to the interesting alternative formulas
\begin{align*}
\nabla_{\vec{q}_i}V(\vec{q}) &= 
   -V_{r_{ij}}(r_{ij},r_{jk},r_{k\ell},r_{ik},r_{j\ell},r_{i\ell})\frac{\vec{r}_{ij}}{r_{ij}} 
   -V_{r_{ik}}(r_{ij},r_{jk},r_{k\ell},r_{ik},r_{j\ell},r_{i\ell})\frac{\vec{r}_{ik}}{r_{ik}} 
   -V_{r_{i\ell}}(r_{ij},r_{jk},r_{k\ell},r_{ik},r_{j\ell},r_{i\ell})\frac{\vec{r}_{i\ell}}{r_{i\ell}}\,,\\[1ex] 
\nabla_{\vec{q}_j}V(\vec{q}) &= 
   \phantom{-} V_{r_{ij}}(r_{ij},r_{jk},r_{k\ell},r_{ik},r_{j\ell},r_{i\ell})\frac{\vec{r}_{ij}}{r_{ij}} 
   -V_{r_{jk}}(r_{ij},r_{jk},r_{k\ell},r_{ik},r_{j\ell},r_{i\ell})\frac{\vec{r}_{jk}}{r_{jk}} 
   -V_{r_{j\ell}}(r_{ij},r_{jk},r_{k\ell},r_{ik},r_{j\ell},r_{i\ell})\frac{\vec{r}_{j\ell}}{r_{j\ell}}\,, \\[1ex] 
\nabla_{\vec{q}_k}V(\vec{q}) &= 
   \phantom{-} V_{r_{ik}}(r_{ij},r_{jk},r_{k\ell},r_{ik},r_{j\ell},r_{i\ell})\frac{\vec{r}_{ik}}{r_{ik}} 
   +V_{r_{jk}}(r_{ij},r_{jk},r_{k\ell},r_{ik},r_{j\ell},r_{i\ell})\frac{\vec{r}_{jk}}{r_{jk}} 
   -V_{r_{k\ell}}(r_{ij},r_{jk},r_{k\ell},r_{ik},r_{j\ell},r_{i\ell})\frac{\vec{r}_{k\ell}}{r_{k\ell}}\,, \\[1ex] 
\nabla_{\vec{q}_\ell}V(\vec{q}) &= 
   \phantom{-} V_{r_{i\ell}}(r_{ij},r_{jk},r_{k\ell},r_{ik},r_{j\ell},r_{i\ell})\frac{\vec{r}_{i\ell}}{r_{i\ell}} 
   +V_{r_{j\ell}}(r_{ij},r_{jk},r_{k\ell},r_{ik},r_{j\ell},r_{i\ell})\frac{\vec{r}_{j \ell}}{r_{j\ell}} 
   +V_{r_{k\ell}}(r_{ij},r_{jk},r_{k\ell},r_{ik},r_{j\ell},r_{i\ell})\frac{\vec{r}_{k\ell}}{r_{k\ell}} \,.
\end{align*}
The representation of the dihedral angle in terms of the distances in lemma~\ref{lem:dihedraldist} might be called the 'distance' definition of the dihedral angle, as a supplement to the cross-product definition and the scalar product definition of the dihedral angle.

Since we have 'distance' definitions of the bond angles as well as the dihedral angles, and since the other potentials are dependent on distances in a natural way, one could evaluate all forces with respect to these potentials in a unified way. We use these representations here in order to use the same technique to construct a discrete gradient. But this unified way to compute the forces, i.e. the gradients of the potentials, might be interesting for the standard Verlet scheme, too.
\subsection{Experiment with butane}

\begin{figure}
\begin{center}
\begin{minipage}{0.49\textwidth}
\renewcommand\lstlistingname{Code fragment}
\begin{lstlisting}[caption={Data butane}, escapechar=\#, label={dgabutane}]
Positions
1 1 0.949003  1.144251 1.100000
2 1 1.000000  1.000000 1.000000
3 1 1.153000  1.000000 1.000000
4 1 1.203997  1.144251 0.900000
5 2 0.375077  0.588336 1.000000
6 2 0.500000  0.500000 1.000000
7 2 0.624923  0.588336 1.000000
8 2 0.749846  0.500000 1.000000
\end{lstlisting}
\end{minipage}
\hspace{-2cm}
\begin{minipage}{0.49\textwidth}
\begin{center}
\hspace*{1cm} \includegraphics[height=2cm]{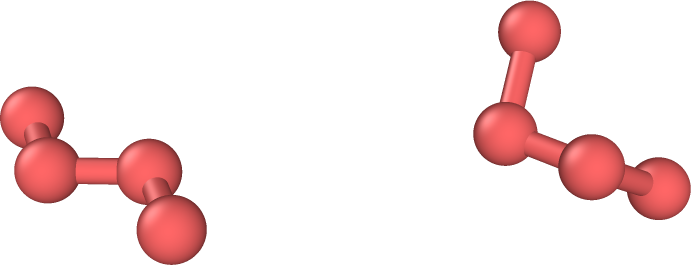}\\[4ex]
\end{center}
\end{minipage}
\end{center}
	\caption{Initial positions for the experiment with two united-atom butane molecules: the positions are given on the left-hand side. The first row numbers the particles. The second row assigns molecule numbers. If this number is the same, a bond is added between the particles. The third to fifth row are the coordinates of the positions. The velocities are set to zero and not given in the code fragment. On the right-hand side, the initial configuration is given as a plot made by the open visualization tool (OVITO), cf. \cite{Stukowski10}.} \label{dgbutaneinitial}
\end{figure}

We use a united-atom model of butane. The Lennard--Jones potential, the bond potential and the angle potential are chosen as before with the parameters given in table~\ref{butanesimdata2}.
And the potential for the torsion potential in terms of the dihedral angle $\phi$ in the IUPAC convention (cf. \cite{IUPAC70}) reads
\begin{align} \label{torsionpotbutane}
  U_t(\phi) &= k_\phi(1.116-1.462\cos \phi-1.578 \cos^2 \phi + 0.368 \cos^3 \phi + 3.156 \cos^4 \phi + 3.788 \cos^5 \phi)\,.
\end{align}
The parameter $k_\phi$ is also given in table~\ref{butanesimdata2} (cf. \cite{Tox87}).
\begin{table}[h]
\begin{center}
{\small
\[
  \begin{array}{lll}
     k_b = 17.5~\frac{\mbox{MJ}}{\mbox{mol}\cdot\mbox{nm}^2}, & r_0=1.53~\mathring{\mbox{A}}, & \mbox{bond potential}, \\[1ex]
     k_\theta = 65~\frac{\mbox{kJ}}{\mbox{mol}}, & \theta_0= 109.47~\mbox{degree}, & \mbox{angle potential}, \\[1ex]
     k_\phi = 8.31451~\frac{\mbox{kJ}}{\mbox{mol}}, & & \mbox{torsion potential}, \\[1ex]
     \sigma = 3.923~\mathring{\mbox{A}}, & \epsilon = 0.5986~\frac{\mbox{kJ}}{\mbox{mol}}, & \mbox{Lennard--Jones potential}.
  \end{array}
\]
}
\end{center}
\caption{Parameters for the simulation of butane} \label{butanesimdata2}
\end{table}
Since we are using more realistic data this time, we first remove the units by scaling. We use
$\tilde{\sigma}=10^{-9}~\mbox{m}=1~\mbox{nm}$, $\tilde{\epsilon}=1~\mbox{kJ}\slash\mbox{mol}$, $\tilde{m}=1~\mbox{u}$ and $\tilde{\alpha}=\tilde{\sigma}\sqrt{\tilde{m}\slash\tilde{\epsilon}} \approx 10^{-12}~\mbox{s}=1~\mbox{ps}$.
The dimensionless quantities then read
\[
 \begin{array}{lllll}
   m'=m\slash \tilde{m},~~~ & \vec{x}_i'=\vec{x}_i\slash \tilde{\sigma},~~~
   & \vec{r}_{ij}'=\vec{r}_{ij}\slash \tilde{\sigma},~~~
   & E'=E\slash \tilde{\epsilon},~~~ & V'=V\slash \tilde{\epsilon}, \\[1ex]
   \sigma'=\sigma\slash \tilde{\sigma},~~~ & \epsilon'=\epsilon\slash \tilde{\epsilon},~~~
   & t'=t\slash \tilde{\alpha}.~~~ & &
 \end{array}
\]
The initial configuration for the atoms in the scaled quantities are shown in figure~\ref{dgbutaneinitial} on the right-hand side.
The molecules are plotted with OVITO (cf. \cite{Stukowski10}).
The first row in code fragment~\ref{dgabutane} on the left-hand side refers to the atom number, the second row to the molecule the atom belongs to.
Two consecutive atoms within the same molecule are connected by a bond. Three consecutive atoms with the same molecule are ruled by the angle potential, and all four atoms of the molecule by the given torsion potential.
That is, this experiment uses all potentials discussed so far. 

\begin{figure}
\begin{center}
\begin{minipage}{0.49\textwidth}
\begin{tikzpicture}

\begin{axis}[
width=8cm, height=6.5cm, xmin=0.5*1e-6, xmax=2e-2, ymin=1e-10, ymax=40, xmode=log, ymode=log, legend pos=south east]
\addplot[mark=*, mark size=1pt, color=red!70!white, thick] coordinates {
( 1.000000000000000021e-02 , 9.785579767440973908e-01 )
( 1.000000000000000021e-03 , 7.562559139512349538e-03 )
( 1.000000000000000048e-04 , 7.573197235626252116e-05 )
( 1.000000000000000082e-05 , 7.573410137896328590e-07 )
( 9.999999999999999547e-07 , 7.568242765514062691e-09 )
};
\addplot[mark=+, mark size=2pt, color=cyan, thick] coordinates {
( 1.000000000000000021e-02 , 2.528082866772254800e+00 )
( 1.000000000000000021e-03 , 1.336888292907852105e-02 )
( 1.000000000000000048e-04 , 1.335418855491838014e-04 )
( 1.000000000000000082e-05 , 1.335396492777906522e-06 )
( 9.999999999999999547e-07 , 1.335720936311215753e-08 )
};
\addplot[mark=diamond*, mark size=2pt, color=orange, thick] coordinates {
( 1.000000000000000021e-02 , 1.395517764947580952e+00 )
( 1.000000000000000021e-03 , 3.444058012604575003e-01 )
( 1.000000000000000048e-04 , 2.967315231564551271e-03 )
( 1.000000000000000082e-05 , 2.941188810934054150e-05 )
( 9.999999999999999547e-07 , 2.763169557573873683e-07 )
};
\addplot[mark=*, mark size=2pt, color=blue!70!white, thick] coordinates {
( 1.000000000000000021e-02 , 1.536506495837130393e+00 )
( 1.000000000000000021e-03 , 1.462735996210861478e+00 )
( 1.000000000000000048e-04 , 1.459200226567393954e-01 )
( 1.000000000000000082e-05 , 1.436120963443324357e-02 )
( 9.999999999999999547e-07 , 1.482232286801451289e-03 )
};
\addplot[color=black, dashed, thick, domain=0.000001:0.01] {10000*x};

\addplot[color=black, thick, domain=0.000001:0.01] {700000*x*x};
\legend{Verlet, midpoint, symmetric DG, simple DG, order1-line, order-2 line};

\end{axis}

\end{tikzpicture}
\end{minipage}
\begin{minipage}{0.49\textwidth}
\input{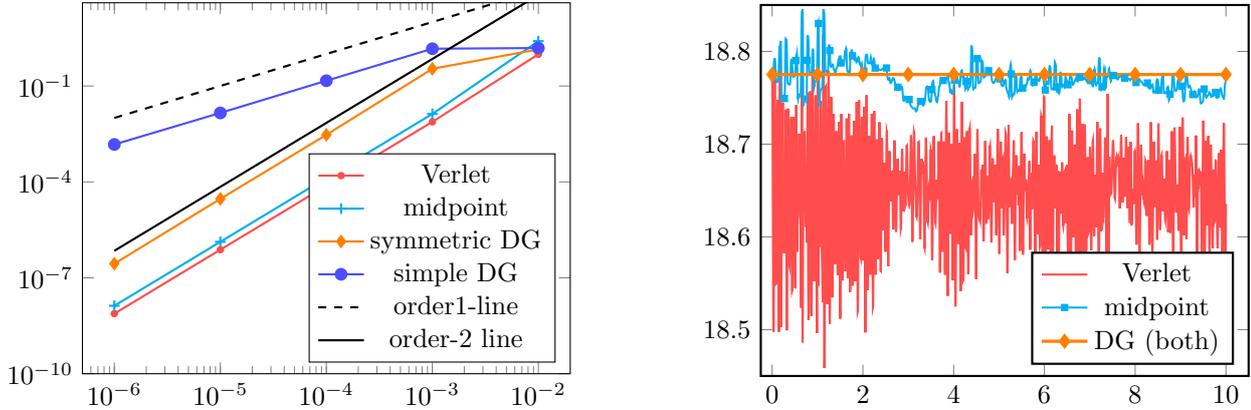}
\end{minipage}
\end{center}
\caption{Results of the experiment with two united-atom butane molecules: the error versus the time step is shown on the left-hand side for the Verlet scheme, the symmetric discrete gradient (DG) scheme, the midpoint rule, and an unsymmetric simple discrete gradient (DG) method. On the right-hand side, the energy is shown over the time span $[0,10]$ for step size $\tau = 0.005$ for the Verlet scheme, the midpoint rule and the two discrete gradient (DG) schemes.} \label{fig:dgdihedral}
\end{figure}

In figure~\ref{fig:dgdihedral} on the left-hand side, one can see that the methods perform as expected with respect to the order. The integration time for the error plot was up to $T=2.0$ with the step sizes indicated on the abscissa. The error measured in the standard Euclidean norm is shown on the ordinate.
If only one of the discrete gradients used in the DG method is unsymmetric then the DG method is of first order. This is shown by the blue circle-marked line.
If all discrete gradients are symmetric, the method is of second order as well as the implicit midpoint rule and the Verlet scheme. For the energy plot on the right-hand side, we computed the solutions up to time $T=10$ with step size $\tau = 0.005$. All DG methods preserve the energy up to round-off error. The implicit midpoint rule and the Verlet scheme deviate from the constant energy at the beginning that should be preserved. 
Since the Verlet scheme deviates significantly more than the implicit midpoint rule, we also calculated the energy with LAMMPS. The energy behaviour turned out to be exactly the same. Setting the NVE ensemble in LAMMPS, the red, solid energy curve is the outcome, which coincides with our own implementation of the Verlet scheme.
The implicit equation in \eqref{velocityDGmethod} is again solved with the simplified Newton method, \eqref{SimpleNewtonmeth}, where the approximate Jacobian $J_F$, cf. \eqref{simpleJac}, includes the full Hessian with respect to the Lennard--Jones potentials and the bond potentials but omits the part with respect to the bond angle and dihedral angle potentials. 
\section{Parallelisation of DG methods} \label{sec:parallel}

In order to show the usefulness of DG methods in molecular dynamics, it is indispensable to take care of parallelisation. Many codes for molecular dynamics simulations are based on the basic parallel treatment of short-range forces (e.g. \cite{Gromacs95,Gromacs01,Lammps95,Lammps22}).

\subsection{Parallelisation of the evaluation of discrete gradients for Lennard--Jones forces with cut-off}
In this section, we discuss the parallelisation of Lennard--Jones forces, induced by the standard potential, cf. \eqref{ljpotential}.
Since we are only interested in the short-range part of the potential and since we will also make use of the Hessian, the Lennard--Jones potential with cut-off function should be twice continuously differentiable. For this reason, we will
use the switching function proposed in \cite{Meietal90} as
\[
 s(r)=
 \left\{
   \begin{array}{ll}
           1, & 0 \leq r < r_{\footnotesize \mathrm{m}},\\
   (1-x)^3(1+3x+6x^2), & r_{\footnotesize \mathrm{m}} \leq r \leq r_{\footnotesize \mathrm{cut}}, \\
           0, & r_{\footnotesize \mathrm{cut}} < r,
   \end{array} 
   \right.  \quad \qquad  x=\frac{r-r_{\footnotesize \mathrm{m}}}{r_{\footnotesize \mathrm{cut}}-r_{\footnotesize \mathrm{m}}}, 
\]
with $r_{\footnotesize \mathrm{m}}=\frac{r_{\footnotesize \mathrm{cut}}}{2}$.
The function and the first two derivatives restricted to $[r_{\footnotesize \mathrm{m}},r_{\footnotesize \mathrm{cut}}]$ read
\[
 \begin{array}{rcl}     
         s(r)   &=& (1-x)^3(1+3x+6x^2) \\
         s'(r)  &=& \frac{1}{r_{\footnotesize \mathrm{cut}}-r_{\footnotesize \mathrm{m}}} 
                  \cdot(-30)\cdot x^2(x-1)^2 \\
         s''(r) &=& \frac{1}{(r_{\footnotesize \mathrm{cut}}-r_{\footnotesize \mathrm{m}})^2}
                  \cdot(-60)\cdot x(x-1)(2x-1)
 \end{array} \, ,\qquad                   
 x=\frac{r-r_{\footnotesize \mathrm{m}}}{r_{\footnotesize \mathrm{cut}}-r_{\footnotesize \mathrm{m}}}\,.
\]
The Lennard--Jones potential $U$ with switching-function, i.e. $V(r)=U(r)\cdot s(r)$, is a twice continuously differentiable short-range potential.
This switching function is available in LAMMPS as {\tt pair\_style lj/mdf}.

\begin{figure}
\begin{center}
\begin{minipage}{0.32\textwidth}
\begin{center}
\begin{tikzpicture}[x=1.0cm,y=1.0cm,scale=1]
   %
	 \filldraw[white!80!black] (1,1) rectangle (4,4);
	 \filldraw[white!60!black,draw=black] (3.75,2.5) arc(0:360:1cm);
         \draw[red, thick] (2.75,2.5) -- (3.457,3.207);
	 \shade [ball color=red, minimum width=0.05cm] (2.75,2.5) circle (3.5pt);
	 \shade [ball color=red, minimum width=0.05cm] (4.25,2.25) circle (3.5pt);
	 \node[red,left] at (2.65,2.5) {{\small $i$}};
	 \node[red,right] at (4.3,2.4) {{\small $j$}};
         \draw[step=1cm,color=black,thick] (0,0) grid (5,5);
         \foreach \x/\y in {0.5/0.5, 4.8/4.8, 1.3/3.6, 3.3/1.2, 
	 1.6/2.9, 2.9/.6, 2.4/2.1, 1.5/4.2, 3.7/3.7, 
	 2.2/3.6, 2.9/1.9, 1.6/1.4, 3.6/4.1, 1.5/3.3, .4/2.0, 2.5/1.4, 
	 4.5/.7, 2.9/4.7}
	 \shade [ball color=blue, minimum width=0.05cm] (\x,\y) circle (3.5pt);
	 %
\end{tikzpicture}
\end{center}
\end{minipage}
\begin{minipage}{0.32\textwidth}	
\begin{center}
\begin{tikzpicture}[x=1cm,y=1cm,scale=1]
	 \filldraw[white!80!black] (1,1) rectangle (4,4);
	 \filldraw[white!60!black,draw=black] (4.15,2.3) arc(0:360:1cm);
         \draw[red, thick] (3.15,2.3) -- (3.857,3.007);
	 \shade [ball color=red] (3.15,2.3) circle (3.5pt);
	 \shade [ball color=red, minimum width=0.05cm] (3.875,2.5) circle (3.5pt);
	 \node[red,above] at (3.25,2.4) {{\small $i'$}};
	 \node[red,below] at (3.65,2.75) {{\footnotesize $j'$}};
   \draw[step=1cm,color=black,thick] (0,0) grid (5,5);
	 \foreach \x/\y in {0.5/0.75, 1.8/0.25, 1.1/3.7, 3.5/1.1, 
	 1.5/2.65, 3.1/.6, 2.3/2.0, 1.2/4, 3.6/3.8, 
	 2.1/3.5, 2.7/1.7, 1.7/1.6, 3.8/3.9, 1.7/3.1, .6/2.3, 2.2/1.2, 
	 4.3/.2, 2.8/4.9}
	 \shade [ball color=blue] (\x,\y) circle (3.5pt);
\end{tikzpicture}
\end{center}
\end{minipage}
\begin{minipage}{0.32\textwidth}
\begin{center}
\begin{tikzpicture}[x=1.0cm,y=1.0cm,scale=1]
	 \filldraw[white!80!black] (0,0) rectangle (5,5);
	 \filldraw[white!60!black,draw=black] (4.75,2.5) arc(0:360:2cm);
         \draw[red, thick] (2.75,2.5) -- (4.16421,3.9142);
	 \shade [ball color=red, minimum width=0.05cm] (2.75,2.5) circle (3.5pt);
	 \shade [ball color=red, minimum width=0.05cm] (4.25,2.25) circle (3.5pt);
	 \node[red,left] at (2.65,2.5) {{\small $i$}};
	 \node[red,right] at (4.3,2.4) {{\small $j$}};
         \draw[step=1cm,color=black,thick] (0,0) grid (5,5);
         \foreach \x/\y in {0.5/0.5, 4.8/4.8, 1.3/3.6, 3.3/1.2, 
	 1.6/2.9, 2.9/.6, 2.4/2.1, 1.5/4.2, 3.7/3.7, 
	 2.2/3.6, 2.9/1.9, 1.6/1.4, 3.6/4.1, 1.5/3.3, .4/2.0, 2.5/1.4, 
	 4.5/.7, 2.9/4.7}
	 \shade [ball color=blue, minimum width=0.05cm] (\x,\y) circle (3.5pt);
\end{tikzpicture}
\end{center}
\end{minipage}
\end{center}
\caption{Linked cell method: Simulation domain is decomposed into square cells of size $r_{\mbox{\tiny cut}}\times r_{\mbox{\tiny cut}}$. The dark-shaded circle is the cut-off radius $r_{\mbox{\tiny cut}}$ about the red particle $i$. 
On the left-hand side, the current time step is shown. For the standard force computation, only the particles in the $3 \times 3$ grid of light-gray cells need to be taken into account. In the middle, the situation at the next time step is shown. The future positions $i'$ of particle $i$ and $j'$ of particle $j$, respectively, are now within the cut-off range. While particle $i$ and $j$ of the current time step do not interact, they interact after they moved to the positions $i'$ and $j'$ in the next time step. In the discrete gradient method, the positions in the next time step are needed to compute the discrete gradient. Therefore, a $5 \times 5$ grid of cells around the cell with the red particle $i$ needs to be taken into account for DG methods, if the particles are assumed to travel not more than $2\slash 3$ of the cut-off in one time step. The cut-off radius of the DG method around the red particle $i$ and the cells that need to be taken into account are shown on the right-hand side.}
\label{linkedcellideaDG}
\end{figure}
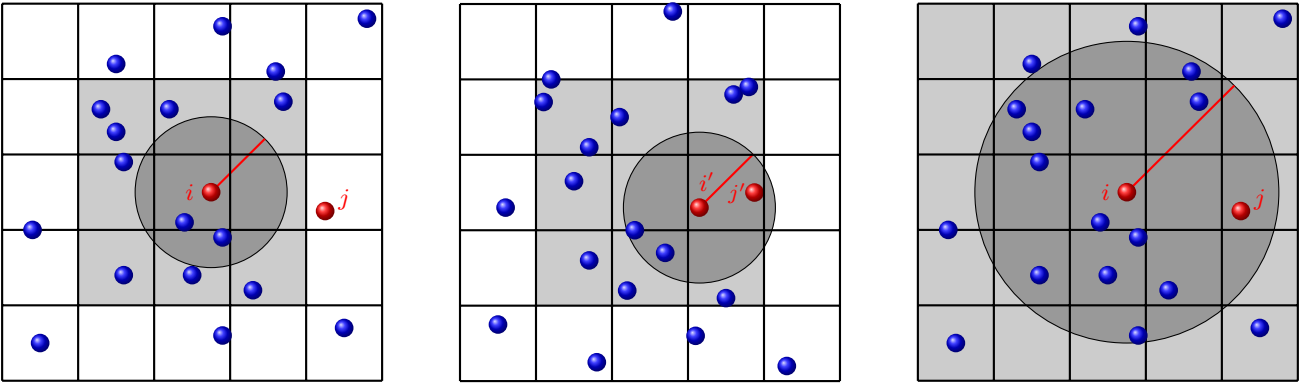

The Lennard--Jones interactions are implemented with the linked cell method as described, e.g., in chapter~3 of \cite{mdgriebel}. Since DG methods are implicit methods, the particle structure is extended by the possible future positions of the particles during the iterative solution of equation \eqref{velocityDGmethod}. Due to this,
Lennard--Jones forces pose a special challenge with DG methods or implicit methods, in general. While a border neighbourhood of one cell is enough for the Verlet scheme, we need a border neighbourhood of two cells for DG methods. The reason is that the particle structure not only carries the position at the given time, but also the position of the next step. Two particles could become close in the next step that are not close in the given time step. This is illustrated in figure~\ref{linkedcellideaDG}. Since possible future positions are needed for the evaluation of the discrete gradient \eqref{dg:pw} in the difference \eqref{dg:pw:diff}, particles with a distance of $2r_\mathrm{\footnotesize cut}$ need to be known, when we assume that the step size is chosen so small that the particles cannot travel for more than $\frac{2}{3}r_\mathrm{\footnotesize cut}$ in a time step.
To sweep all particles in a border neighbourhood of two cells of a given cell is enough to catch these events. As an alternative, one might use cells with dimensions larger than $2r_\mathrm{\footnotesize cut}$. Then a border neighbourhood with width one cell would be enough for the larger cells. We will stick with the smaller cells. There are some computations, for example the potential with the cut-off function, where we only need a border neighbourhood of one cell. Hence, using the smaller cells saves a bit of computing time.

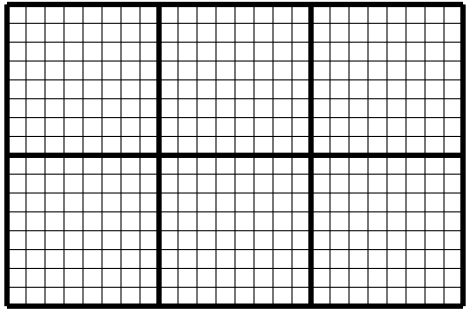
\begin{figure}
\begin{center}
\begin{tikzpicture}[scale=0.25]
\foreach \x in {0,...,24}
  \draw (\x,0) -- (\x,16);
\foreach \x in {0,8,16,24}
  \draw[line width=2pt] (\x,0) -- (\x,16);
\foreach \y in {0,...,16}
  \draw (0,\y) -- (24,\y);
\foreach \y in {0,8,16}
  \draw[line width=2pt] (0,\y) -- (24,\y);
\end{tikzpicture}
\end{center}
\caption{Decomposition of the simulation domain: the domain $\Omega$ is divided in cells indicated by thin lines. The domain is subdivided further into six subdomains as indicated by the thick lines. Each of the six subdomains is assigned to its own processor that handles the computations in this subdomain according to the linked cell method.} \label{subdomdecompDG}
\end{figure}

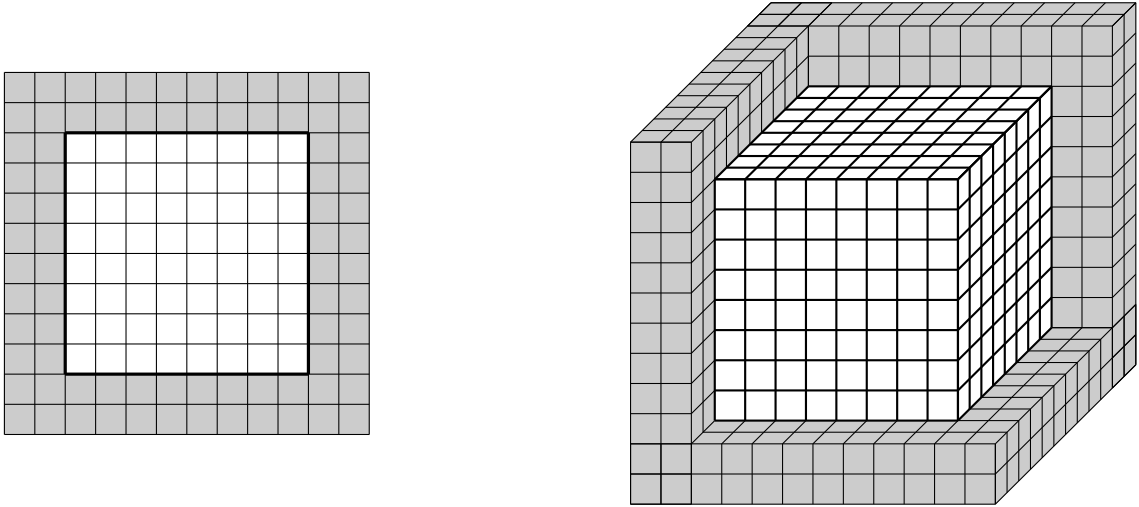
\begin{figure}
\begin{center}
\begin{minipage}{0.49\textwidth}
\begin{center}
\begin{tikzpicture}[scale=0.4]
\draw[fill, white!80!black] (0,0) rectangle (2,12);
\draw[fill, white!80!black] (10,0) rectangle (12,12);
\draw[fill, white!80!black] (0,0) rectangle (12,2);
\draw[fill, white!80!black] (0,10) rectangle (12,12);
\foreach \x in {0,...,12}
  \draw (\x,0) -- (\x,12);
\foreach \x in {2,10}
  \draw[very thick] (\x,2) -- (\x,10);
\foreach \y in {0,...,12}
  \draw (0,\y) -- (12,\y);
\foreach \y in {2,10}
  \draw[very thick] (2,\y) -- (10,\y);
\end{tikzpicture}
\end{center}
\end{minipage}
\begin{minipage}{0.49\textwidth}
\begin{center}
\begin{tikzpicture}[scale=0.4]
\draw[fill, white!80!black] (0,0,12) rectangle (2,12,12);
\draw[fill, white!80!black] (0,12,12) -- (2,12,12) -- (2,12,0) -- (0,12,0) -- (0,12,12);
\draw[fill, white!80!black] (0,0,12) -- (12,0,12) -- (12,2,12) -- (0,2,12) -- (0,0,12);
\draw[fill, white!80!black] (12,0,12) -- (12,0,0) -- (12,2,0) -- (12,2,12) -- (12,0,12);
\draw[fill, white!80!black] (12,0,0) -- (12,12,0) -- (12,12,2) -- (12,0,2) -- (12,0,0);
\draw[fill, white!80!black] (12,12,0) -- (0,12,0) -- (0,12,2) -- (12,12,2) -- (12,12,0);
\draw[fill, white!80!black] (2,2,12) -- (12,2,12) -- (12,2,10) -- (2,2,10) -- (2,2,12);
\draw[fill, white!80!black] (12,2,12) -- (10,2,12) -- (10,2,2) -- (12,2,2) -- (12,2,12);
\draw[fill, white!80!black] (12,2,2) -- (12,2,2) -- (12,12,2) -- (10,12,2) -- (10,2,2);
\draw[fill, white!80!black] (12,10,2) -- (12,12,2) -- (2,12,2) -- (2,10,2) -- (12,10,2);
\draw[fill, white!80!black] (2,12,2) -- (2,12,12) -- (2,10,12) -- (2,10,2) -- (2,12,2);
\draw[fill, white!80!black] (2,12,12) -- (2,12,10) -- (2,2,10) -- (2,2,12) -- (2,12,12);
\draw (0,12,0) -- (0,12,12);
\draw (1,12,0) -- (1,12,12);
\draw (2,12,0) -- (2,12,12);
\foreach \x in {0,...,12}
  \draw (0,12,\x) -- (2,12,\x);
\draw (0,12,12) -- (0,0,12);
\draw (1,12,12) -- (1,0,12);
\draw (2,12,12) -- (2,0,12);
\foreach \x in {0,...,12}
  \draw (0,\x,12) -- (2,\x,12);
\draw (0,0,12) -- (12,0,12);
\draw (0,1,12) -- (12,1,12);
\draw (0,2,12) -- (12,2,12);
\foreach \x in {0,...,12}
  \draw (\x,0,12) -- (\x,2,12);
\draw (12,0,12) -- (12,0,0);
\draw (12,1,12) -- (12,1,0);
\draw (12,2,12) -- (12,2,0);
\foreach \x in {0,...,12}
  \draw (12,0,\x) -- (12,2,\x);
\draw (12,0,0) -- (12,12,0);
\draw (12,0,1) -- (12,12,1);
\draw (12,0,2) -- (12,12,2);
\foreach \x in {0,...,12}
  \draw (12,\x,0) -- (12,\x,2);
\draw (0,12,0) -- (12,12,0);
\draw (0,12,1) -- (12,12,1);
\draw (0,12,2) -- (12,12,2);
\foreach \x in {0,...,12}
  \draw (\x,12,0) -- (\x,12,2);
\foreach \x in {2,...,10}
  \draw (2,10,\x) -- (2,12,\x);
\draw (2,11,2) -- (2,11,12);
\foreach \x in {2,...,10}
  \draw (2,\x,10) -- (2,\x,12);
\draw (2,11,2) -- (12,11,2);
\foreach \x in {3,...,10}
  \draw (\x,2,10) -- (\x,2,12);
\draw (11,2,12) -- (11,2,2);
\foreach \x in {2,...,10}
  \draw (10,2,\x) -- (12,2,\x);
\draw (11,2,2) -- (11,12,2);
\foreach \x in {2,...,10}
  \draw (10,\x,2) -- (12,\x,2);
\draw (12,2,11) -- (2,2,11);
\foreach \x in {2,...,10}
  \draw (\x,10,2) -- (\x,12,2);
\draw (2,2,11) -- (2,12,11);
\foreach \x in {2,...,10}
  \foreach \y in {10}
    \draw[thick] (\x,\y,2) -- (\x,\y,10);
\foreach \x in {10}
  \foreach \y in {2,...,10}
    \draw[thick] (\x,\y,2) -- (\x,\y,10);
\foreach \x in {2,...,10}
  \foreach \y in {10}
    \draw[thick] (\x,2,\y) -- (\x,10,\y);
\foreach \x in {10}
  \foreach \y in {2,...,10}
    \draw[thick] (\x,2,\y) -- (\x,10,\y);
\foreach \x in {2,...,10}
  \foreach \y in {10}
    \draw[thick] (2,\x,\y) -- (10,\x,\y);
\foreach \x in {10}
  \foreach \y in {2,...,10}
    \draw[thick] (2,\x,\y) -- (10,\x,\y);
\end{tikzpicture}
\end{center}
\end{minipage}
\end{center}
\caption{Subdomain and border neighbourhood: On the left-hand side, a typical subdomain is shown as the square with the white cells. The processor responsible for this domain needs the data from the adjacent domains. A border neighbourhood of cells, given in gray, is added to the subdomain. The border neighbourhood is filled with copies of the particles in the adjacent domains. The adjacent processors have to send the data to the current processor by messages, according to the message passing paradigm standardised in the message passing interface (MPI), cf. \cite{MPI}. For implicit methods and DG methods, a border neighbourhood of two cells is necessary in order to compute the forces by the linked cell method. On the right-hand side, the situation for a three-dimensional subdomain is illustrated.} \label{subdomneighimplpar}
\end{figure}

For the parallelisation, domain decomposition is used. If we decompose the two-dimensional domain in figure~\ref{subdomdecompDG} in six larger parts based on the cells given by the linked cell method, every processor only knows the particles in its domain. Since each processor might run on its own node with its own memory in a cluster computer, the access to the data of adjacent processors is not immediate (cf. \cite{clustercomp}). The standard technique is to extend the domain that the processor is handling by further cells, the border neighbourhood, and to retrieve the necessary information from the adjacent processors in these cells. The current processor also has to send the information needed by the adjacent processors from his domain to the neighbours. 
For the exchange of the data between processors, message passing is used. That is, the processors send messages to each other. This is standardised in the message passing interface (MPI) (cf. \cite{MPI}).
Due to the discussion above, we need to extend the processor's domain by a border neighbourhood of the size of two cells as shown in figure~\ref{subdomneighimplpar} on the left hand side. We will conduct three-dimensional experiments. The corresponding domain and its neighbourhood is shown in figure~\ref{subdomneighimplpar} on the right-hand side.

\subsection{Collision of two bodies}
\begin{figure} 
\input{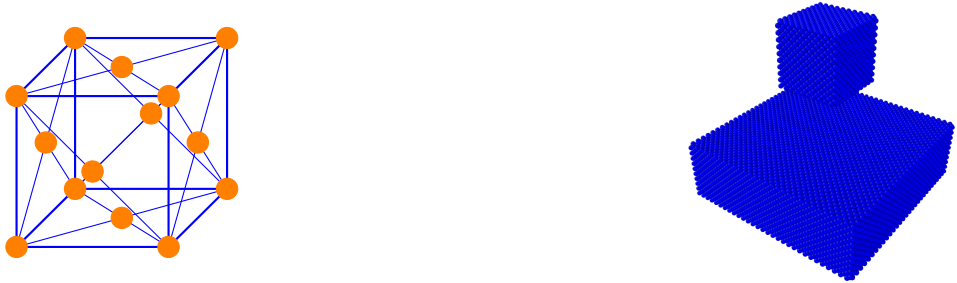}
\caption{Sketch of a face-centered cubic lattice on the left-hand side. On the right-hand side, the initial positions of the particles that form the smaller and the large body are shown in a plot made by the open visualization tool (OVITO), cf. \cite{Stukowski10}.} \label{fcclatticepic}
\end{figure}

\begin{figure}
\begin{center}
\begin{minipage}{0.49\textwidth}
\begin{center}
\includegraphics[width=6cm, height=6cm]{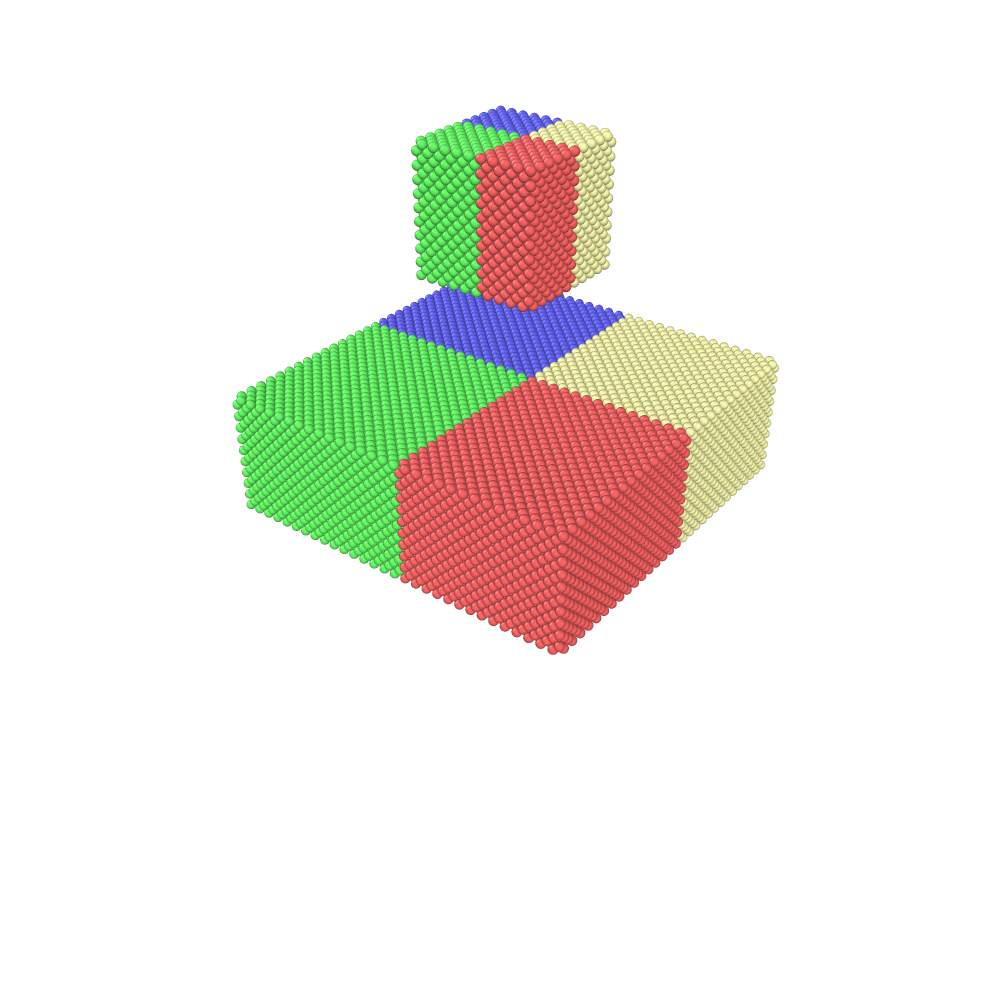}\\[4ex]
\end{center}
\end{minipage}
\begin{minipage}{0.49\textwidth}
\begin{center}
\includegraphics[width=6cm, height=6cm]{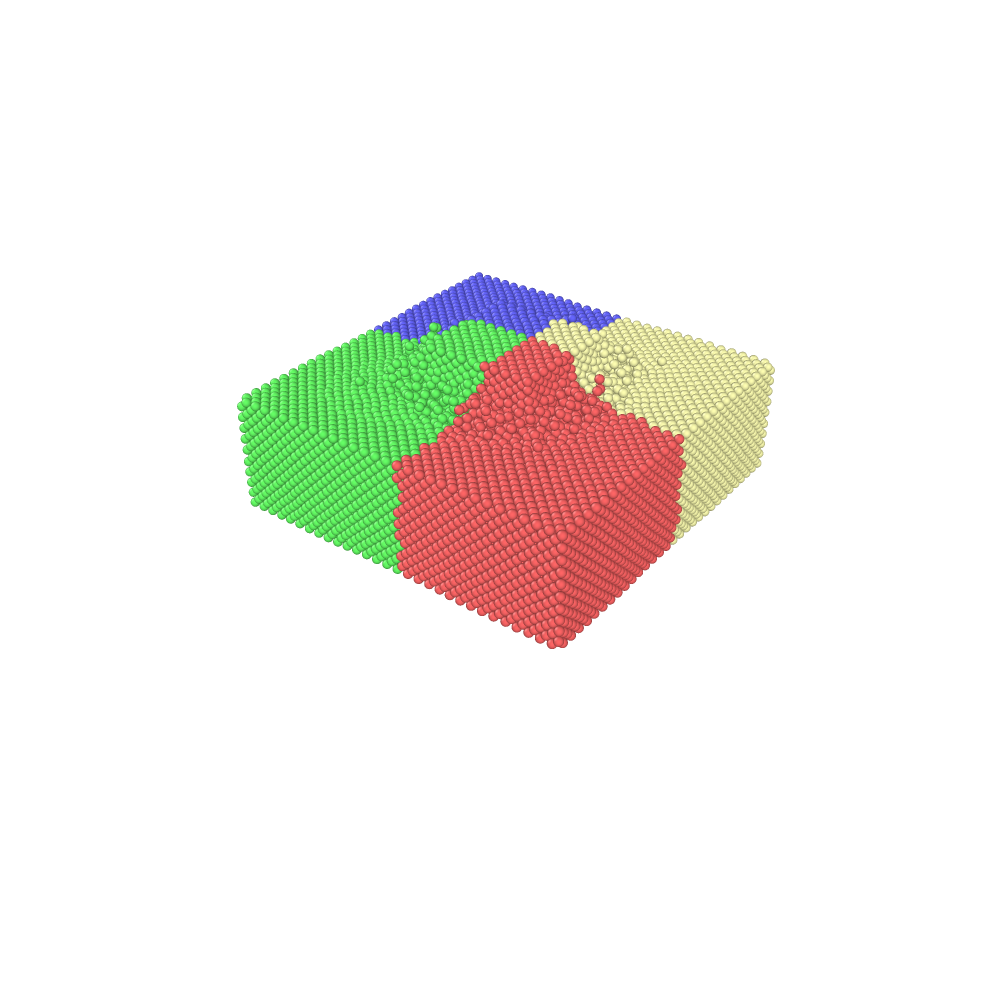}\\[4ex]
\end{center}
\end{minipage}
\vspace*{-1.25cm}

\begin{minipage}{0.49\textwidth}
\begin{center}
\includegraphics[width=6cm, height=6cm]{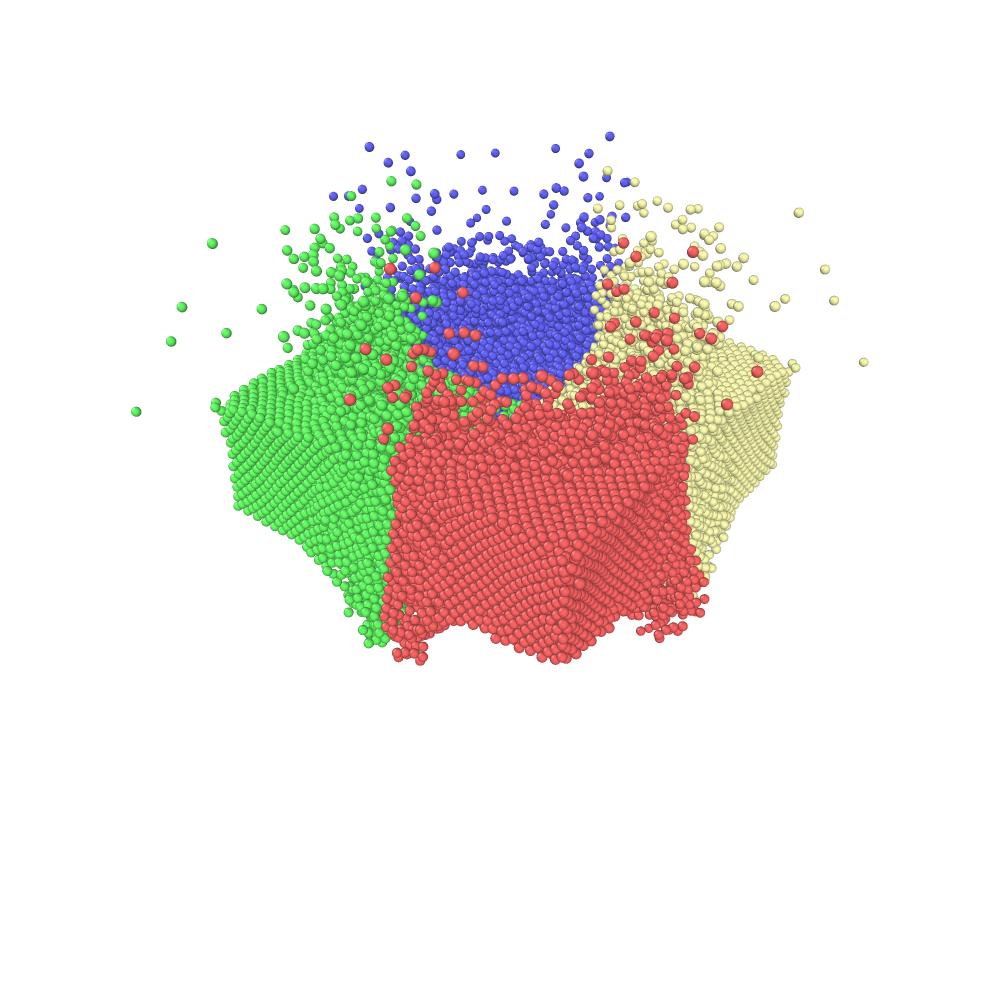}\\[4ex]
\end{center}
\end{minipage}
\begin{minipage}{0.49\textwidth}
\begin{center}
\includegraphics[width=6cm, height=6cm]{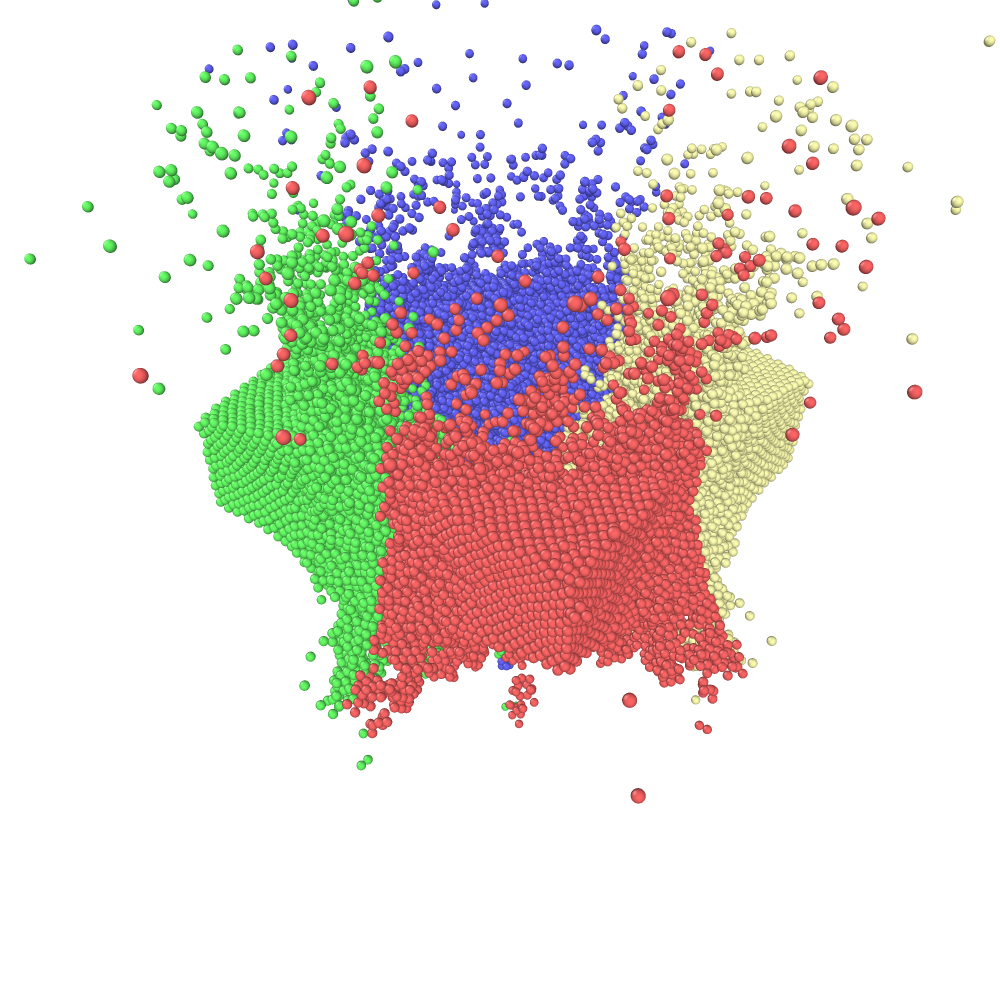}\\[4ex]
\end{center}
\end{minipage}

\vspace*{-1cm}

\begin{minipage}{0.49\textwidth}
\begin{center}
\includegraphics[width=6cm, height=6cm]{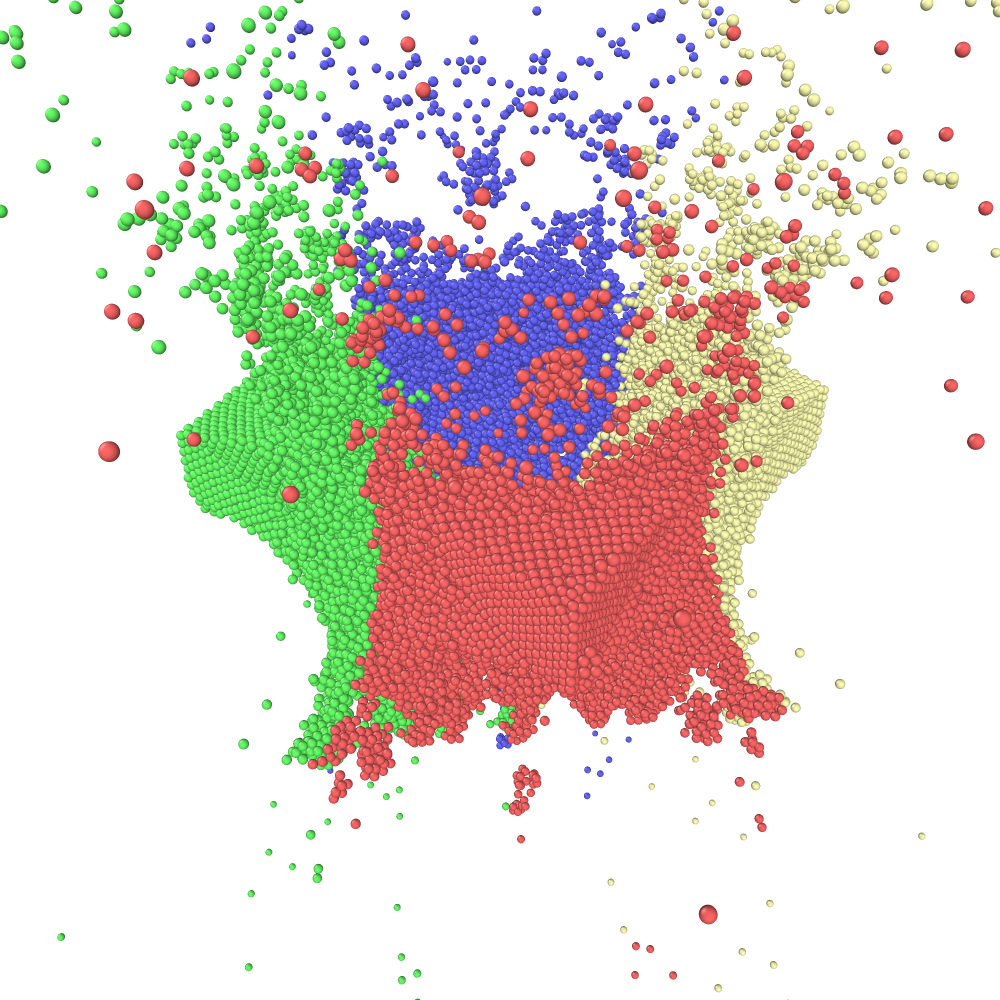}\\[4ex]
\end{center}
\end{minipage}
\begin{minipage}{0.49\textwidth}
\begin{center}
\includegraphics[width=6cm, height=6cm]{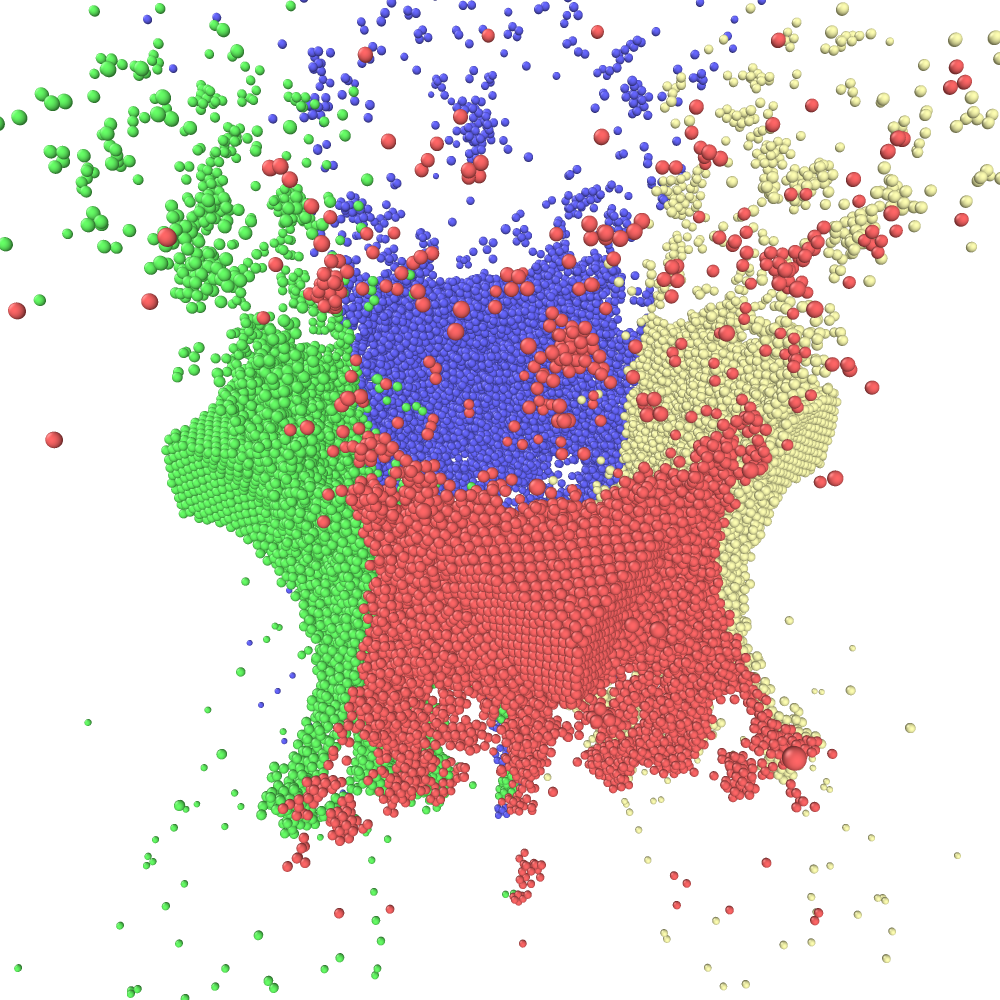}\\[4ex]
\end{center}
\end{minipage}
\end{center}
\caption{Simulation of the collision of two bodies at times $t=0,4,8,12,16,20$ from top left to bottom right. The simulation is run with four processors identified by color. The same color means that the same processor is handling the particles.} \label{crash4proc}
\end{figure}

As a test problem for the parallelisation, we use the collision of two bodies as described in section 4.5.1 
of \cite{mdgriebel}. In figure~\ref{fcclatticepic} on the right-hand side, the initial configuration of the experiment is shown. The two bodies consist of $10 \times 10 \times 10$ and $10 \times 30 \times 30$ cubic cells, with $4000$ and $36000$ particles, respectively. Each cell contains four particles in a face-centered cubic grid,
 as shown on the left-hand side of figure~\ref{fcclatticepic}. The four particles, which are counted for this cell, are the lower left corner as well as the particles in the center of the front, bottom and left face. The others belong to regular grids spanned by these four particles. The shortest distance in the grid is $2^{1\slash 6}\sigma$, according to the equilibrium length of the Lennard--Jones potential.
At the beginning of the simulation, the smaller body moves with a high velocity $\vec{v}$ towards the resting larger body. The simulation cell of size $[0,150\sigma]^3$ is equipped with periodic boundary conditions. All data are given in table~\ref{crash3dsimdata}.
\begin{table}[h]
\begin{center}
{\small
\[
 \begin{array}{rclp{0.5cm}rclp{0.5cm}rcl}
   L_1 &=& 150\sigma, && L_2 &=& 150\sigma, & &L_3 &=& 150\sigma \\[1.5ex]
         \varepsilon &=& 120, && \sigma &=& 3.4, &&&& \\[1.5ex]
           m &=& 39.95, && \vec{v} &=& (0,0,-20.4), &&&& \\[1.5ex]
        \mbox{dist. part.} &=& 2^{1\slash 6}\sigma, && N_1 &=& 4000, && N_2 &=& 36000 \\[1.5ex]
                r_{\mbox{\scriptsize cut}}&=& 2.5\sigma, && \tau&=& 0.001 &&&&
  \end{array}
\]
}
\end{center}
\caption{Parameters for the simulation of the collision} \label{crash3dsimdata}
\end{table}
Several snapshots of the simulation are shown in figure~\ref{crash4proc}. The simulation has been conducted with four processors. For our proof-of-concept implementation, we used four standard personal computers (PCs) as processors that had been connected by a one gigabit local area network (1 GB LAN). The particles in figure~\ref{crash4proc} are color-coded with respect to the processor that handles the particles.
\begin{figure}
\begin{center}
\input{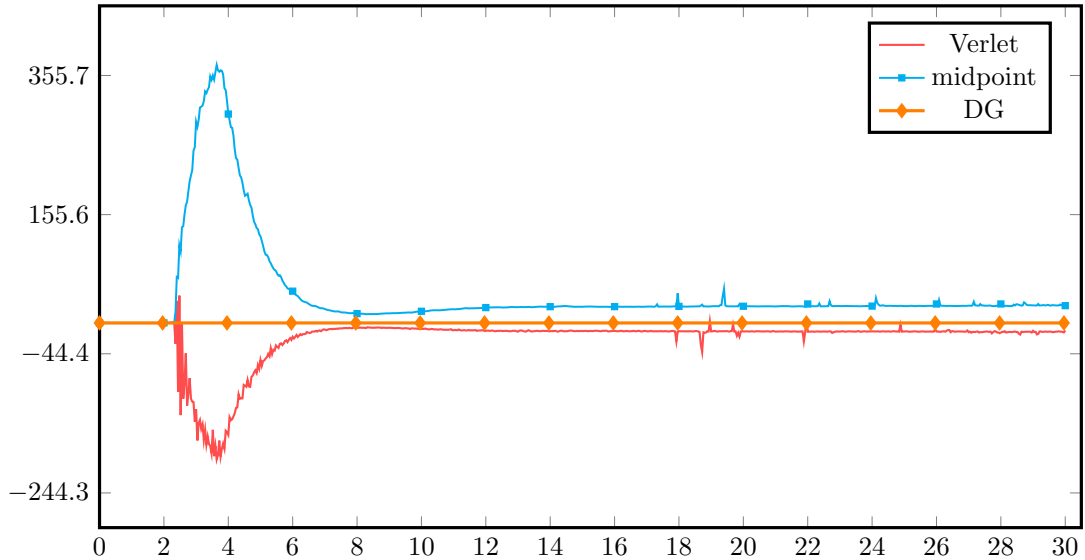}
\end{center}
\caption{Energy preservation for the collision of two bodies: the deviation from the exact energy is shown over the time span $[0,30]$ for step size $\tau=0.001$ for the Verlet scheme, the midpoint rule and the discrete gradient (DG) scheme. All three methods have been computed in parallel with four processors on a cluster computer, cf. \cite{clustercomp}.} \label{dgcrashenergy}
\end{figure}
The simulation is run up to (scaled) time $T=30$. All results are given in scaled quantities. In figure~\ref{dgcrashenergy}, one can see that the DG method preserves the energy very well. The implicit midpoint rule overestimates the energy when the small body penetrates the larger body. The Verlet scheme underestimates the energy during this phase. If the larger body is destroyed, the energy computed by the midpoint rule decreases and the energy computed by the Verlet scheme increases. In order to check our computation, we also conduct the same experiment with LAMMPS and only one processor (no parallelisation). The observed energy behaviour of the independent Verlet implementation in LAMMPS shows exactly the same energy curve, including the small peaks later on. That is, only the discrete gradient method is able to simulate a true microcanonical ensemble. The implicit equation in \eqref{velocityDGmethod} is solved with the Newton method, \eqref{Newtonmeth}, with the full Jacobian \eqref{NewtonmethJac}. The action of the Jacobian on a vector is directly computed. This way, one can transfer the linked cell method to the computation of the action of the Jacobian on a vector. The total linear momentum is preserved by all three methods. The total angular momentum can not be preserved, by any of the methods, due to the periodic boundary conditions. Only for free space or repelling boundary conditions, methods can preserve the total angular momentum in a MD simulation.
\subsection{Parallelisation of the evaluation of discrete gradients for bonded forces}
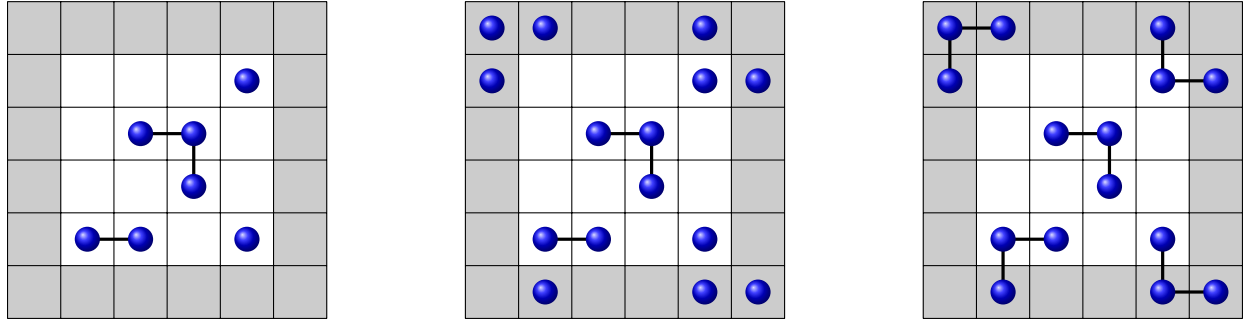
\begin{figure}
\begin{center}	
\begin{center}
\begin{minipage}{0.32\textwidth}
\begin{center}
\begin{tikzpicture}[scale=0.7]
	 \filldraw[white!80!black] (0,0) rectangle (6,1);
	 \filldraw[white!80!black] (0,1) rectangle (1,6);
	 \filldraw[white!80!black] (5,1) rectangle (6,6);
	 \filldraw[white!80!black] (1,5) rectangle (5,6);
        \draw[step=1cm,color=black] (0,0) grid (6,6);
	 \draw[very thick] (1.5,1.5) -- (2.5,1.5);
	 \node[circle,shading=ball,minimum width=0.2cm] (ball) at (1.5,1.5) {};
	 \node[circle,shading=ball,minimum width=0.2cm] (ball) at (2.5,1.5) {};
	 \node[circle,shading=ball,minimum width=0.2cm] (ball) at (4.5,1.5) {};
	 \draw[very thick] (2.5,3.5) -- (3.5,3.5);
	 \draw[very thick] (3.5,3.5) -- (3.5,2.5);
	 \node[circle,shading=ball,minimum width=0.2cm] (ball) at (2.5,3.5) {};
	 \node[circle,shading=ball,minimum width=0.2cm] (ball) at (3.5,3.5) {};
	 \node[circle,shading=ball,minimum width=0.2cm] (ball) at (3.5,2.5) {};
	 \node[circle,shading=ball,minimum width=0.2cm] (ball) at (4.5,4.5) {};
\end{tikzpicture}
\end{center}
\end{minipage}
\begin{minipage}{0.32\textwidth}
\begin{center}
\begin{tikzpicture}[scale=0.7]
	 \filldraw[white!80!black] (0,0) rectangle (6,1);
	 \filldraw[white!80!black] (0,1) rectangle (1,6);
	 \filldraw[white!80!black] (5,1) rectangle (6,6);
	 \filldraw[white!80!black] (1,5) rectangle (5,6);
   \draw[step=1cm,color=black] (0,0) grid (6,6);
	 \draw[very thick] (1.5,1.5) -- (2.5,1.5);
	 \node[circle,shading=ball,minimum width=0.2cm] (ball) at (1.5,0.5) {};
	 \node[circle,shading=ball,minimum width=0.2cm] (ball) at (1.5,1.5) {};
	 \node[circle,shading=ball,minimum width=0.2cm] (ball) at (2.5,1.5) {};
	 \node[circle,shading=ball,minimum width=0.2cm] (ball) at (4.5,1.5) {};
	 \node[circle,shading=ball,minimum width=0.2cm] (ball) at (4.5,0.5) {};
	 \node[circle,shading=ball,minimum width=0.2cm] (ball) at (5.5,0.5) {};
	 \draw[very thick] (2.5,3.5) -- (3.5,3.5);
	 \draw[very thick] (3.5,3.5) -- (3.5,2.5);
	 \node[circle,shading=ball,minimum width=0.2cm] (ball) at (2.5,3.5) {};
	 \node[circle,shading=ball,minimum width=0.2cm] (ball) at (3.5,3.5) {};
	 \node[circle,shading=ball,minimum width=0.2cm] (ball) at (3.5,2.5) {};
	 \node[circle,shading=ball,minimum width=0.2cm] (ball) at (0.5,4.5) {};
	 \node[circle,shading=ball,minimum width=0.2cm] (ball) at (0.5,5.5) {};
	 \node[circle,shading=ball,minimum width=0.2cm] (ball) at (1.5,5.5) {};
	 \node[circle,shading=ball,minimum width=0.2cm] (ball) at (4.5,5.5) {};
	 \node[circle,shading=ball,minimum width=0.2cm] (ball) at (4.5,4.5) {};
	 \node[circle,shading=ball,minimum width=0.2cm] (ball) at (5.5,4.5) {};
\end{tikzpicture}
\end{center}
\end{minipage}
\begin{minipage}{0.32\textwidth}
\begin{center}
\begin{tikzpicture}[scale=0.7]
	 \filldraw[white!80!black] (0,0) rectangle (6,1);
	 \filldraw[white!80!black] (0,1) rectangle (1,6);
	 \filldraw[white!80!black] (5,1) rectangle (6,6);
	 \filldraw[white!80!black] (1,5) rectangle (5,6);
   \draw[step=1cm,color=black] (0,0) grid (6,6);
	 \draw[very thick] (1.5,0.5) -- (1.5,1.5);
	 \draw[very thick] (1.5,1.5) -- (2.5,1.5);
	 \node[circle,shading=ball,minimum width=0.2cm] (ball) at (1.5,0.5) {};
	 \node[circle,shading=ball,minimum width=0.2cm] (ball) at (1.5,1.5) {};
	 \node[circle,shading=ball,minimum width=0.2cm] (ball) at (2.5,1.5) {};
	 \draw[very thick] (4.5,1.5) -- (4.5,0.5);
	 \draw[very thick] (4.5,0.5) -- (5.5,0.5);
	 \node[circle,shading=ball,minimum width=0.2cm] (ball) at (4.5,1.5) {};
	 \node[circle,shading=ball,minimum width=0.2cm] (ball) at (4.5,0.5) {};
	 \node[circle,shading=ball,minimum width=0.2cm] (ball) at (5.5,0.5) {};
	 \draw[very thick] (2.5,3.5) -- (3.5,3.5);
	 \draw[very thick] (3.5,3.5) -- (3.5,2.5);
	 \node[circle,shading=ball,minimum width=0.2cm] (ball) at (2.5,3.5) {};
	 \node[circle,shading=ball,minimum width=0.2cm] (ball) at (3.5,3.5) {};
	 \node[circle,shading=ball,minimum width=0.2cm] (ball) at (3.5,2.5) {};
	 \draw[very thick] (0.5,4.5) -- (0.5,5.5);
	 \draw[very thick] (0.5,5.5) -- (1.5,5.5);
	 \node[circle,shading=ball,minimum width=0.2cm] (ball) at (0.5,4.5) {};
	 \node[circle,shading=ball,minimum width=0.2cm] (ball) at (0.5,5.5) {};
	 \node[circle,shading=ball,minimum width=0.2cm] (ball) at (1.5,5.5) {};
	 \draw[very thick] (4.5,5.5) -- (4.5,4.5);
	 \draw[very thick] (4.5,4.5) -- (5.5,4.5);
	 \node[circle,shading=ball,minimum width=0.2cm] (ball) at (4.5,5.5) {};
	 \node[circle,shading=ball,minimum width=0.2cm] (ball) at (4.5,4.5) {};
	 \node[circle,shading=ball,minimum width=0.2cm] (ball) at (5.5,4.5) {};
\end{tikzpicture}
\end{center}
\end{minipage}
\end{center}
\end{center}
\caption{Reconnection and separation of atoms in the border neighbourhood: If the subdomain of the processor receives atoms that are bonded within a molecule, the processor needs to reattach the atoms at the correct sites. This is illustrated here from left to right in a two-dimensional example. First the subdomain with the knowledge of the processor is shown. In the middle, the situation is shown that results after the processor has received the particles necessary for the computation of the bonded forces from the adjacent processors. With the help of the molecule and particle numbers, the processor reattaches the molecules in the correct way. The result is illustrated on the right-hand side. After the computation of the bonded forces, the process is reversed from right to left.} \label{cellbondedideaDG}
\end{figure}

The parallelisation of the bonded forces is simpler in the sense that it works in the standard way, which is illustrated in figure~\ref{cellbondedideaDG}. On the left-hand side of figure~\ref{cellbondedideaDG}, the domain with its border neighbourhood is shown. Then the adjacent processors send the necessary particles, while this processor sends the particles needed by neighbouring processors. This is the state in the middle of figure~\ref{cellbondedideaDG}.
Then the particles in the border neighbourhood are linked based on the molecule numbers and atom numbers stored with the particles. After that, the current processor sees all data needed to compute the discrete gradients with respect to bonded forces for the particles in its domain. Recall, that the new positions are also stored with the corresponding particles. This state is illustrated on the right-hand side of figure~\ref{cellbondedideaDG}. After their computation, the discrete gradients are also stored with the particles. Then, the particles in the border neighbourhood are separated again, which brings us back to the state in the middle of the figure. Finally, the particles in the border neighbourhood are deleted, because they are no longer needed. This procedure also works for the matrix-free computation of the Hessian of the bonded potentials times a vector.
\subsection{Experiment with butane}
\begin{figure}
\begin{center}
\begin{minipage}{0.49\textwidth}
\begin{center}
\includegraphics[width=8.5cm, height=8.5cm]{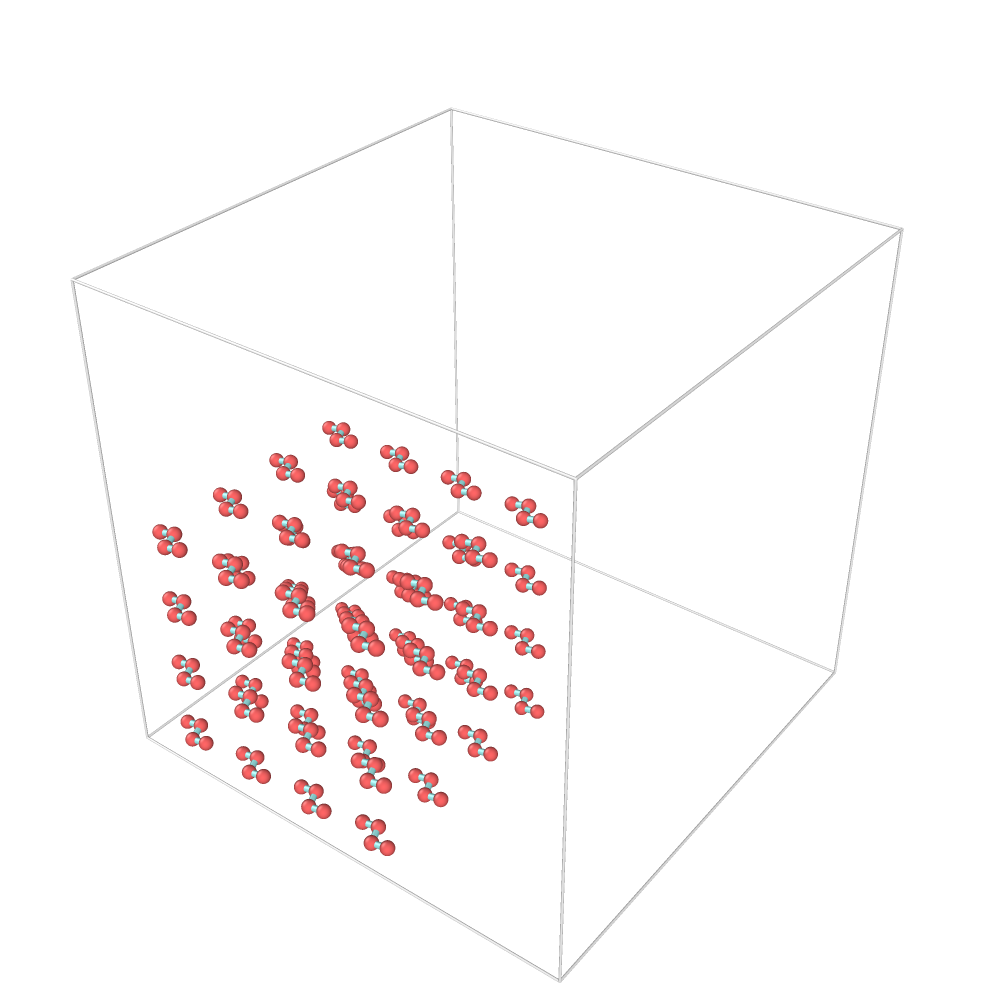}\\[4ex]
\end{center}
\end{minipage}
\begin{minipage}{0.49\textwidth}
\begin{center}
\includegraphics[width=8.5cm, height=8.5cm]{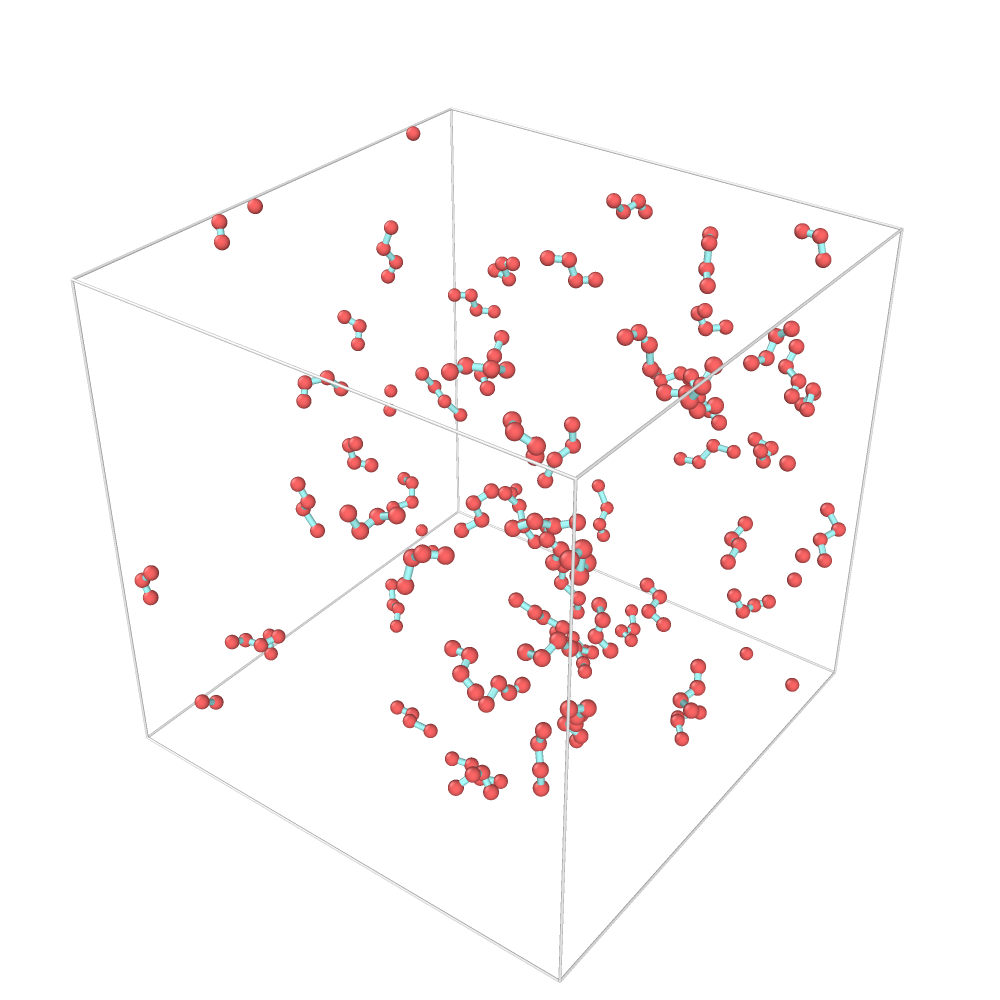}\\[4ex]
\end{center}
\end{minipage}
\end{center}
	\caption{Simulation of butane in a periodic box at times $t=0,20$. The plot is made by the open visualization tool (OVITO), cf. \cite{Stukowski10}} \label{butaneboxproc}
\end{figure}

\begin{figure}
\begin{center}
\includegraphics[width=14.5cm]{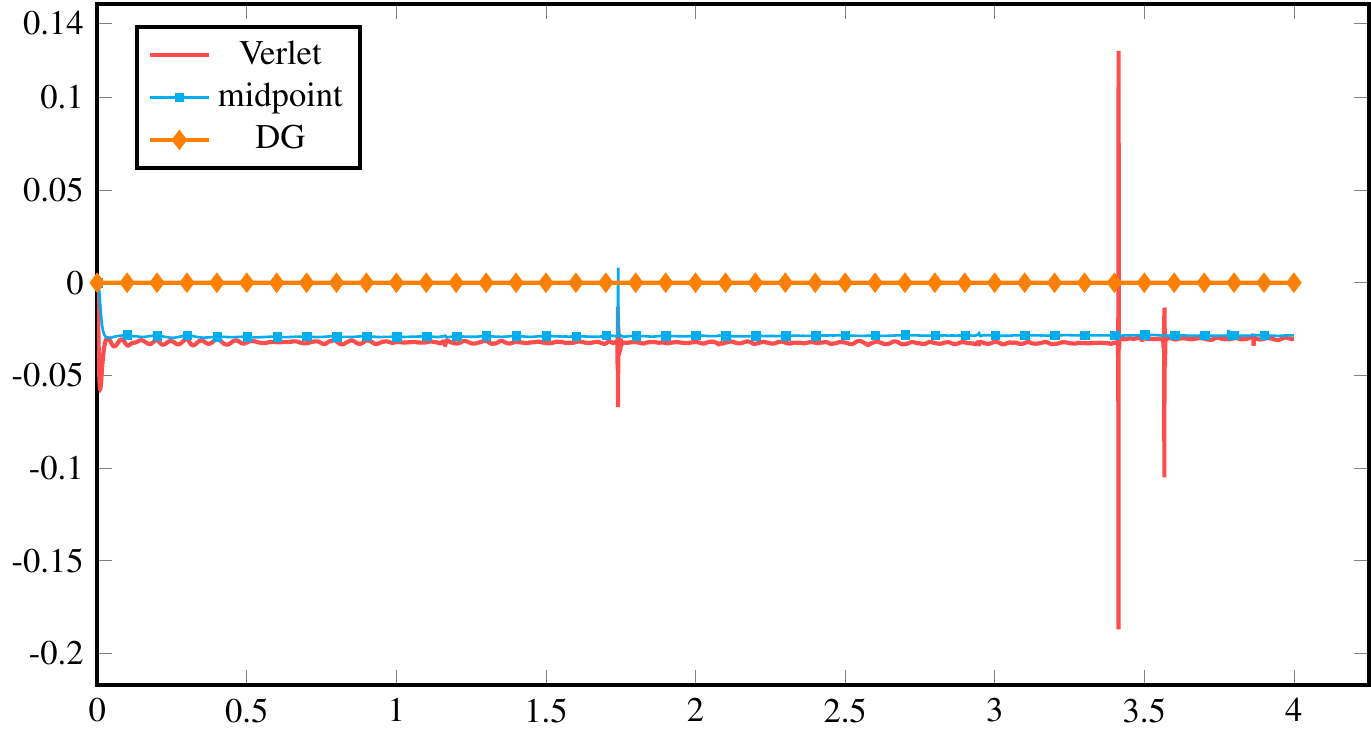}\\[4ex]
\end{center}
\caption{Energy preservation during the simulation of $64$ united-atom butane molecules: the deviation from the exact energy is shown over the time span $[0,4]$ with step size $\tau=0.0001$ for the Verlet scheme, the midpoint rule, and the discrete gradient (DG) method.} \label{dgbutaneenergy}
\end{figure}

We run an experiment with $64$ united-atom butane molecules. The initial condition and a snapshot of the simulation at time $T=4$ can be seen in figure~\ref{butaneboxproc}. The data for the potentials have been chosen as in table~\ref{butanesimdata2}.	Also the potentials are chosen as before, i.e., the angle potential is given in equation~\eqref{anglepotcossq} and the torsion potential is given in equation~\eqref{torsionpotbutane}. The simulation is set in a periodic box of size $4r_{\footnotesize \mathrm{cut}}$, $r_{\footnotesize \mathrm{cut}}=2.5\sigma$, and $\sigma$ as in table~\ref{butanesimdata2} in scaled variables. The whole simulation is scaled as before in the butane simulation.
The deviation from the exact energy over the simulation time with step size $\tau = 0.0001$ is shown in figure~\ref{dgbutaneenergy}. 
The simulation has been run with four processors. 
While the DG methods preserve the energy up to round-off error, the implicit midpoint scheme and the Verlet scheme deviate from the constant energy. The peaks in the energy of the Verlet scheme are real. We also computed the energy for the given initial value with LAMMPS with one processor. This simulation reproduced exactly the same peaks as our code. This means that the particles do not evolve with respect to a genuine NVE ensemble at these peaks. The solid red line is the closest a standard molecular dynamics package with the Verlet scheme and setting the ensemble to the NVE ensemble can get.   
If one wishes more accuracy with respect to energy preservation, discrete gradient methods are an interesting alternative.
\section{Conclusion} \label{sec:conclusion}
This work shows that all standard short-range interactions in a classical conservative molecular dynamics simulation can be computed by discrete gradient methods. 
These methods reliably preserve the total energy in the system, along with the total linear momentum and the total angular momentum in free space simulations. 
The simple and unified idea to construct the discrete gradients is to express all standard short-range interactions in terms of distances between atoms.
The new discrete gradients for the dihedral angle potentials also suggest an interesting way to compute the gradient of dihedral angle potentials based on distances for the use in standard time integration schemes. Furthermore, the discrete gradient methods can be parallelised. We proposed the necessary changes to the linked cell method for the parallel evaluation of the discrete gradients with respect to truncated Lennard--Jones potentials as well as the necessary changes for bonded forces. As a result, the proposed DG methods can be computed in parallel.    
\section*{Acknowledgments} The authors would like to thank the Isaac Newton Institute for Mathematical Sciences, Cambridge, for support and
hospitality during the programme ``Geometry, compatibility and structure preservation in computational differential equations'' where work on this paper was undertaken. This work was supported by EPSRC grant no EP/K032208/1. In addition, this work was supported by the German Science Foundation (DFG) under project GRK 2450.

\end{document}